\documentclass{amsart}

\usepackage{float}                       % Float placement
\usepackage{amsmath, amsthm, amssymb}      % AMS math packages
\usepackage{graphicx}                    % For inserting images
\usepackage{longtable,todonotes}         % Long tables and to-do notes
\usepackage{tabularx, colortbl}           % Advanced table features
\usepackage{caption}                     % Customize captions
\usepackage{mathtools, mathrsfs, cancel}   % Extended math tools
\usepackage{stmaryrd}                     % Extra symbol sets
\usepackage[inline]{enumitem}             % Inline lists

\usepackage{etex}                         % Extended registers
\usepackage{etoolbox}                     % Enhanced macro tools

\usepackage{fancyhdr}                     % Fancy headers/footers
\usepackage{hyperref}                     % Hyperlinks in document
\usepackage[capitalize, nameinlink, noabbrev]{cleveref}  % Smart cross-references

\PassOptionsToPackage{table,dvipsnames,svgnames}{xcolor}
\usepackage{xcolor}                       % Extended color support

\usepackage{tcolorbox}                    % Colored boxes
\tcbuselibrary{breakable, skins, theorems}  % tcolorbox libraries
\usepackage[final]{pdfpages}              % Include external PDF pages

\usepackage{tikz, tikz-cd, tikz-qtree}      % TikZ and diagram drawing
\usepackage{tkz-euclide}                   % Euclidean geometry drawings
\usetikzlibrary{%
    arrows, backgrounds, decorations, decorations.markings, 
    decorations.pathreplacing, decorations.pathmorphing, external, 
    calc, hobby, intersections, knots, math, matrix, patterns, 
    positioning, through, shapes%
}
\tikzset{>=latex}
\pgfdeclarelayer{background}
\pgfdeclarelayer{foreground}
\pgfsetlayers{background,main,foreground}
\usepackage{subcaption} % Required for \subcaptionbox
\usetikzlibrary{tqft}
\usetikzlibrary{arrows.meta, positioning}
\newtheorem{thm}{Theorem}[section]

\newtheorem{lem}[thm]{Lemma}

\newtheorem{prop}[thm]{Proposition}

\theoremstyle{definition}
\newtheorem{remark}[thm]{Remark}

\newtheorem{defi}[thm]{Definition}

\def\liesl{\mathfrak{sl}}          % Special linear group Lie algebra
\DeclareMathOperator{\ORTH}{O}      % Orthogonal group
\DeclareMathOperator{\PSL}{PSL}     % Projective special linear group
\DeclareMathOperator{\SL}{SL}       % Special linear group
\DeclareMathOperator{\SO}{SO}       % Special orthogonal group
\usepackage{dsfont}  % For \mathds (blackboard bold)
\DeclareMathOperator{\tr}{tr}

\DeclareMathOperator{\Isom}{Isom}

\makeatletter
\renewenvironment{proof}[1][\proofname]{%
   \par\pushQED{\qed}\normalfont%
   \topsep6\p@\@plus6\p@\relax%
   \trivlist\item[\hskip\labelsep\bfseries#1\@addpunct{.}]%
   \ignorespaces%
}{%
   \popQED\endtrivlist\@endpefalse%
}
\makeatother

\DeclareMathOperator{\osys}{osys}

\author{Nhat Minh Doan}
\address{Department of Mathematics, National University of Singapore, Singapore,\newline
\& Institute of Mathematics, Vietnam Academy of Science and Technology, Vietnam}

\email{minh.dn@nus.edu.sg, dnminh@math.ac.vn}
\urladdr{https://sites.google.com/view/dnminh}

\author{Khanh Le}
\address{Department of Mathematics, Rice University, USA}
\email{khanh.le.math@gmail.com}
\urladdr{https://sites.google.com/view/khanhqle}



\title{Some arithmetic aspects of ortho-integral surfaces}

\begin{document}

\maketitle
\begin{abstract}
We investigate ortho-integral (OI) hyperbolic surfaces with totally geodesic boundaries, defined by the property that every orthogeodesic (i.e. a geodesic arc meeting the boundary perpendicularly at both endpoints) has an integer cosh-length. We prove that while only finitely many OI surfaces exist for any fixed topology, infinitely many commensurability classes arise as the topology varies. Moreover, we completely classify OI pairs of pants and OI one-holed tori, and show that their doubles are arithmetic surfaces of genus 2 derived from quaternion algebras over $\mathbb{Q}$. 
\end{abstract}

\begin{center}
%\tableofcontents
\end{center}

\section{Introduction}

For a hyperbolic $n$-manifold $X$ with non-empty totally geodesic boundary, we consider the set of \emph{orthogeodesics} which are geodesic arcs meeting $\partial X$ perpendicularly at both endpoints. The \emph{ortho-length spectrum}, $\mathcal{O}(X)$, is the ordered multi-set of length of orthogeodesics. The spectrum was introduced by Basmajian in the study of embedded totally geodesic hypersurfaces in hyperbolic manifolds \cite{basmajian1993orthogonal} which relates the ortho-length spectrum and the geometry of the manifold. A consequence of Basmajian's work in dimension two is that the total boundary length of a hyperbolic surface $X$ with geodesic boundary can be recovered from the ortho-length spectrum, in particular:
\begin{equation}
\label{eq:BasmajianIdentity}
\ell(\partial X) = \sum_{L \in \mathcal{O}(X) } 2 \log(\coth(L/2)).    
\end{equation}
Since the work of Basmajian, the ortho‑length spectrum has been studied extensively in relation to geometric identities and hyperbolic geometry of the surface. In \cite{bridgeman2011orthospectra}, Bridgeman established an identity connecting the volume of the unit tangent bundle, the ortho‑length spectrum, and the dilogarithm function, with Calegari later offering an alternative proof \cite{calegari2011bridgeman}. Bridgeman and Kahn then generalized this result by showing that the volume of a hyperbolic manifold with nonempty totally geodesic boundary can be recovered from its ortho‑length spectrum \cite{bridgeman2010hyperbolic}. In \cite{bridgeman2016identities}, Bridgeman and Tan described the common framework underlying these and other geometric identities. For related well-known identities and applications, see \cite{mcshane1998simple,mirzakhani2007simple,tan2008generalized,bridgeman2014moments,luo2014a,doan2023measuring,parlier2020geodesic} and the references therein. Recently, Basmajian, Parlier, and Tan \cite{basmajian2020prime} generalized Basmajian’s identity to hyperbolic surfaces with cusps, geodesic boundaries, and cone points by considering prime orthogeodesics in the concave core obtained by cutting off natural collars from the surface.

In connection with the hyperbolic geometry of a surface, the ortho-length spectrum has been studied as an analog of the conventional length spectrum. Masai and McShane \cite{masai2023systoles} investigated the rigidity and flexibility of the ortho-length spectrum, showing that it uniquely determines the hyperbolic structure on a one-holed torus, constructing examples of non-isometric hyperbolic surfaces with identical ortho-length spectra, and proving that for any fixed hyperbolic surface with non-empty geodesic boundaries, only finitely many non-isometric surfaces share the same ortho-length spectrum. For further recent progress in this direction, see \cite{quellec2024orthospectrum}. Bavard \cite{bavard2005anneaux} proved that for a hyperbolic surface $X$ with genus $g$ and $n$ boundary components, the length of the shortest orthogeodesic and the total length of the boundary of $X$ cannot both be too large:
\begin{equation}
\label{eq:BavaradInequality}
2\sinh\left(\frac{\osys(X)}{2}\right)\sinh\left(\frac{\ell(\partial X)}{12|\chi(X)|}\right)\leq 1    
\end{equation}
where $\osys(X)$ is the length of the shortest orthogeodesic on $X$. Expanding on \cite{bavard2005anneaux}, Basmajian and Fanoni studied the local maxima of the orthosystole function $\osys(X)$ on the moduli space of hyperbolic surfaces with prescribed boundary lengths \cite{basmajian2024orthosystoles}. 

The starting point of this article is the work of the first author in \cite{doan2021}, where he studied the combinatorial structure of the set of orthogeodesics. In particular, he showed that for a compact hyperbolic surface $X$ with nonempty totally geodesic boundary, the set of orthogeodesics on $X$ can be parametrized by a finite collection of rooted planar trivalent trees which depends on a decomposition of $X$ into right-angled hexagons. Using this tree structure, he gave a recursive method to compute the set of hyperbolic cosine lengths of all orthogeodesics, with multiplicity, called the \emph{ortho cosh-length spectrum}:  
\[
\mathcal O_{\cosh}(X)
:=
\bigl\{\cosh\bigl(\ell(\mu)\bigr)\;\bigm|\;\mu\text{ is an orthogeodesic on }X\bigr\},
\] thus producing a combinatorial proof of Basmajian’s identity analogous to Bowditch’s combinatorial proof of the McShane identity. When the surface is a pair of pants or a one-holed torus, there always exists a hexagon decomposition consisting of two isometric right-angled hexagons. Due to the symmetry in this decomposition, the cosh-lengths of orthogeodesics on a pair of pants or a one-holed torus can be realized as solutions to a single non-homogeneous quadratic equation in three variables.

This connection motivates our investigation into the arithmetic properties of the ortho cosh-length spectrum of hyperbolic surfaces. Since the hyperbolic cosine of half the translation length of a hyperbolic isometry in $\PSL_2(\mathbb{R})$ equals half the absolute value of the trace of its matrix, the cosh-length of an orthogeodesic is an analog of the trace of the hyperbolic isometry translating along a closed geodesic. We consider the following definition:

\begin{defi}[\cite{doan2021}]
\label{def:OI}
A hyperbolic surface \(X\) is called \emph{ortho integral (OI)} if $\mathcal{O}_{\cosh}(X) \subset  \mathbb{Z}$, that is, the cosh of the length of every orthogeodesic on \(X\) is an integer.
\end{defi}

In that work, the first author provided examples of OI pairs of pants and OI one-holed tori and posed the problem of classifying OI surfaces. In this paper, we present some arithmetic characterizations and classification results for these surfaces. Our first main result is a finiteness theorem:

\begin{thm}\label{thm:FinitenessOfOISurfaces}
The set of OI surfaces with fixed genus \(g\) and \(n>0\) boundary components is nonempty and finite.
\end{thm}

The nonemptiness in \cref{thm:FinitenessOfOISurfaces} is established constructively in \cref{prop:InfinitelyManyOISurfaces} by gluing certain identical OI pairs of pants together without twist along their boundaries, either by gluing a pair of pants with itself or with another, to form an OI surface of the desired signature \((g,n)\).
 The finiteness is proved via a compactness argument in the moduli space of hyperbolic surfaces with nonempty boundary. In fact, we establish this result for a slightly broader class of hyperbolic surfaces (see \cref{thm:FinitenessOfdAOISurfaces} for the full statement). By providing explicit lower bounds for the systole and bounds for the lengths of the boundary components of OI surfaces, we obtain the desired finiteness. A key step in this argument is given by Bavard's inequality (\cref{eq:BavaradInequality}) together with the following integrality property of closed geodesics on OI surfaces:

\begin{thm}\label{thm:intergerTrace}
Every closed geodesic \(\mu\) on an OI surface satisfies
\[
\tr_{\SL}^2(\mu) \in \mathbb{Z},
\]
where \(|\tr_{\SL}(\mu)| = 2\cosh\bigl(\ell(\mu)/2\bigr).\)
\end{thm}

This integrality extends to a larger class of hyperbolic surfaces (see \cref{rem:IntegralityLargerClass}). We proved it using the recursive formula in \cite{doan2021} for computing the cosh-lengths of certain orthogeodesics on a pair of pants and an immersion of the pair of pants into the surface. Moreover, through this integrality property and a more careful analysis, we obtain a complete classification of OI pairs of pants and OI one-holed tori (see \cref{thm:classificationOIpants} and \cref{thm:classificationOItorus} for the full statements), confirming that all such surfaces were already found in \cite{doan2021}.

Another promising approach to proving \cref{thm:FinitenessOfOISurfaces} is to show that the doubles of these surfaces are arithmetic. In fact, a consequence of the volume formula for a maximal order defining an arithmetic Fuchsian group is that there exist only finitely many arithmetic surfaces with bounded area (see, e.g., \cite{borel1981commensurability} and \cite[Chapter 11]{maclachlan2003book}). Therefore, if every OI surface doubles to yield an arithmetic surface, finiteness of OI surfaces would follow immediately. Although we did not pursue this strategy directly, we investigated the arithmetic properties of OI pairs of pants and OI one-holed tori and verified that the doubles of these OI surfaces are indeed arithmetic surfaces of genus 2 derived from quaternion algebras over \(\mathbb{Q}\) (see \cref{prop:Double,tab:OIPantsInvariants,tab:OIToriInvariants}). It would be interesting to further explore the relationship between arithmetic surfaces and the doubling construction of OI surfaces.

Recall that two hyperbolic surfaces are in the same commensurability class if they share an isometric finite sheeted covering. Motivated by the question of whether there exist infinitely many commensurability classes of pseudomodular surfaces (that is, surfaces of finite area that are not commensurable with the modular surface and that have cusp sets exactly equal to \(\mathbb{Q}\cup\{\infty\}\); see \cite{long2002pseudomodular} and \cite{lou2018hyperbolic} for the answer to this question), the first author asked whether there exist infinitely many commensurability classes of OI surfaces (\cite{doan2021}). In \cref{sec:ConstructingOIsurfaces}, we answer this question affirmatively by constructing infinitely many pairwise non-commensurable OI surfaces obtained by gluing together certain OI pairs of pants (see \cref{fig:GluingPants} for an example).

\begin{thm}\label{thm:InfinitelyManyCommensurabilityClasses}
    There are infinitely many commensurability classes of OI surfaces.
\end{thm}
While most non‑commensurability proofs invoke arithmetic or trace‑field invariants (e.g. via quaternion algebras or length‑spectrum comparisons), our argument is entirely geometric.  Since the surfaces have nonempty totally geodesic boundary, each gluing of the fixed pair of pants strictly enlarges the convex core and raises its maximal inscribed‑disk radius, producing a strictly increasing sequence of inradii and hence an infinite family of pairwise non‑commensurable surfaces.

\subsection{Organization}
In \cref{sec:Prelim}, we review preliminary results on Mumford's compactness criterion, the orthospectrum, hyperbolic identities, and the arithmetic invariants of Fuchsian groups. In \cref{sec:Finiteness}, we prove the integrality property of closed geodesics on OI surfaces and establish the finiteness of the set of OI surfaces with a fixed topological type. In \cref{sec:Classification of OI pants and OI one-holed tori}, we classify OI pairs of pants and OI one-holed tori. Finally, \cref{sec:Ortho-integral surfaces and quadratic forms} explores connections between OI surfaces and quadratic forms by verifying the arithmeticity of the doubles of OI pairs of pants and one-holed tori, proving the nonemptiness of the set of OI surfaces with a fixed topological type, and establishing the infiniteness of commensurability classes of OI surfaces.

\subsection{Acknowledgement}
We thank Hugo Parlier, Alan Reid, and Ser Peow Tan for their invaluable discussions and support. We are especially grateful to Greg McShane for his suggestions and for providing the initial data and the idea of \cref{prop:PantsAndOneHoledTorusFromRootFlipping}, and to Trong Toan Dao for his assistance with the proof of \cref{thm:classificationOIpants}. We thank the anonymous referee for several valuable comments and suggestions that have improved the paper. We also thank the organizers of the Computational Aspects of Thin Groups workshop at the Institute for Mathematical Sciences - National University of Singapore, for hosting the event where this project began. The first-named author is funded by the Vietnam National Foundation for Science and Technology Development (NAFOSTED) under grant number 101.04-2023.33 and by the Singapore National Research Foundation (NRF) under grant E-146-00-0029-01. The second-named author is funded by the AMS Travel Grant.

\section{Preliminary}
\label{sec:Prelim}
\subsection{Notations and terminologies}
Throughout this paper, the term \emph{surface} always refers to a two-dimensional, orientable topological manifold with a negative Euler characteristic. We denote by $\Sigma_{g,n}$ a surface of genus $g$ with $n$ boundary components. Given a hyperbolic structure on a surface $\Sigma_{g,n}$, we get a discrete subgroup of $\PSL_2(\mathbb{R})$ that is the image of the holonomy representation of $\pi_1(\Sigma_{g,n})$ that corresponds to the hyperbolic structure. When $n \geq 1$, $\pi_1(\Sigma_{g,n})$ is free, and the representation $\rho$ can be lifted to a representation also call $\rho: \pi_1(\Sigma_{g,n})\to\SL_2(\mathbb{R}).$ For a primitive closed geodesic \(\gamma\) on a hyperbolic surface with length \(\ell(\gamma)\), we denote by \(\tr_{\SL}(\gamma)\) the trace of its corresponding hyperbolic element in the image of $\rho$, we have\begin{equation}\label{eq:Length&Trace}
    |\tr_{\SL}(\gamma)| = 2\cosh \frac{\ell(\gamma)}{2}.
    \end{equation}

\subsection{Moduli space of surfaces, systole, and Mumford's compactness criterion} 

The moduli space \(\mathcal{M}_{g,n}\) is defined as the space of complete hyperbolic structures with geodesic boundaries on the surface \(\Sigma_{g,n}\), considered up to isometries that fix each boundary component individually. To construct \(\mathcal{M}_{g,n}\), we begin with the Teichmüller space \(\mathcal{T}(\Sigma_{g,n})\), given by
\[
\mathcal{T}(\Sigma_{g,n}) = \left\{ (X, f) \,\middle|\,
\begin{array}{l}
X \text{ is a hyperbolic surface with geodesic boundaries}, \\
f: \Sigma_{g,n} \rightarrow X \text{ is a homeomorphism}
\end{array}
\right\} \Big/ \sim,
\]
where two pairs \((X, f)\) and \((X', f')\) are equivalent, \((X, f) \sim (X', f')\), if there exists an isometry \(m: X \rightarrow X'\) such that the composition \((f')^{-1} \circ m \circ f\) is isotopic to the identity on \(\Sigma_{g,n}\). Here, the homeomorphism \(f\) is called the \emph{marking}.

The mapping class group \(\mathrm{Mod}(\Sigma_{g,n})\) consists of orientation-preserving diffeomorphisms of \(\Sigma_{g,n}\) that fix each boundary component individually. This group acts on \(\mathcal{T}(\Sigma_{g,n})\) by changing the marking:
\[
\phi \cdot [(X, f)] = [(X, f \circ \phi^{-1})], \quad \text{for } \phi \in \mathrm{Mod}(\Sigma_{g,n}).
\]
This action is properly discontinuous, and the quotient
\[
\mathcal{M}_{g,n} = \mathcal{T}(\Sigma_{g,n}) / \mathrm{Mod}(\Sigma_{g,n})
\]
yields the moduli space. Thus, \(\mathcal{M}_{g,n}\) parametrizes the complete hyperbolic structures with geodesic boundaries on \(\Sigma_{g,n}\), up to isometries that fix the boundary components individually.

The \emph{systole}, \(\operatorname{sys}(X)\), of a hyperbolic surface \(X\) is the length of its shortest non-boundary closed geodesic. Systoles play a fundamental role in the study of hyperbolic surfaces and their moduli spaces. For instance, the original statement of Mumford's compactness criterion \cite{mumford1971remark} characterizes compact subsets of the moduli space of closed surfaces based on systolic lower bounds. The following theorem extends this criterion to surfaces with geodesic boundary components.

\begin{thm}[Extended version of Mumford's Compactness Criterion]
\label{thm:MumfordCompactnesCriterion}
A subset \(K \subset \mathcal{M}_{g,n}\) is compact if and only if there exist constants \(\epsilon > 0\) and \(M > 0\) such that for every surface \(X \in K\), the following conditions hold:
\[
  \operatorname{sys}(X) \;\ge\; \epsilon 
  \quad\text{and}\quad
  \epsilon \;\le\; \ell_X(\delta_i) \;\le\; M
  \quad\text{for all geodesic boundary components } \delta_i.
\]
\end{thm}

% \begin{remark}

%    In \cite[Lemma~6.3]{basmajian2024orthosystoles}, a sequence of surfaces is constructed with orthosystoles (the lengths of the shortest orthogeodesics) uniformly bounded below and systoles approaching zero, showing that an analogue of Mumford’s compactness criterion using orthosystole is impossible.

% \end{remark}

% \begin{defi}[Moduli Space]
% The \textbf{moduli space} $\mathcal{M}(\Sigma_{g,n,p})$ is the quotient of the Teichmüller space by the action of the mapping class group:
% \[
% \mathcal{M}(\Sigma_{g,n,p}) = \frac{\mathcal{T}(\Sigma_{g,n,p})}{\mathrm{Mod}(\Sigma_{g,n,p})}.
% \]
% \end{defi}

\subsection{Orthospectrum of hyperbolic surfaces}

An \emph{orthobasis} (also known as a \emph{hexagon decomposition} \cite{basmajian2024orthosystoles} or \emph{truncated triangulation} \cite{ushijima1999canonical}) $\mathcal{H}$ of a surface \( \Sigma_{g,n} \) is a maximal collection of pairwise disjoint, embedded, essential, simple arcs connecting boundary components, possibly the same boundary component. Here an arc is essential if it is not homotopic rel endpoints into one of the boundary components of the surface. 

In a hyperbolic structure \( X \in \mathcal{T}(\Sigma_{g,n}) \), the orthogeodesic representatives of \( \mathcal{H} \) decompose \( X \) into right-angled hexagons. If every three orthogeodesics at the boundary of each hexagon are pairwise distinct, then \( \mathcal{H} \) is a \emph{standard orthobasis}. We use the term orthobasis as it serves as input for computing the \emph{cosh-length spectrum}, which consists of the hyperbolic cosines of orthogeodesic lengths \cite{doan2021}. By Euler characteristic considerations, an orthobasis of a surface $X$ with signature \( (g, n) \) contains \( 3|\chi(X)| = 6g - 6 + 3n \) arcs and divides the surface into \( 2|\chi(X)| = 4g - 4 + 2n \) topological disks. Ushijima \cite{ushijima1999canonical} showed that:

\begin{thm}\label{thm:HexDecompositionHomeoTeich}
Let \( \mathcal{H} = \{\alpha_1, \dots, \alpha_{6g - 6 + 3n}\} \) be an orthobasis of a surface \( \Sigma_{g,n} \). Then the map
\[
\varphi_{\mathcal{H}} : \mathcal{T}(\Sigma_{g,n}) \to \mathbb{R}_{>0}^{6g - 6 + 3n}, \quad X \mapsto \bigl( \ell_X(\alpha_i) \bigr)_{i = 1}^{6g - 6 + 3n},
\]
is a homeomorphism.
\end{thm}
\noindent
An alternative proof for \cref{thm:HexDecompositionHomeoTeich} can be found also in \cite[Theorem 3.1]{basmajian2024orthosystoles}.

\begin{figure}
    \begin{subfigure}{0.25\textwidth}
    \centering
    \includegraphics[{width = \linewidth}]{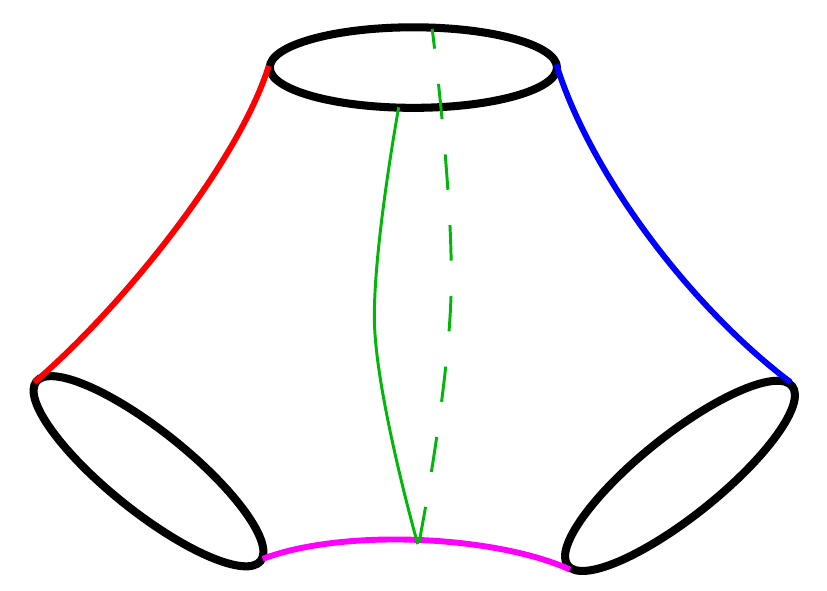}        
    \end{subfigure}
    ~
    \begin{subfigure}{0.25\textwidth}
    \centering
    \includegraphics[{width = \linewidth}]{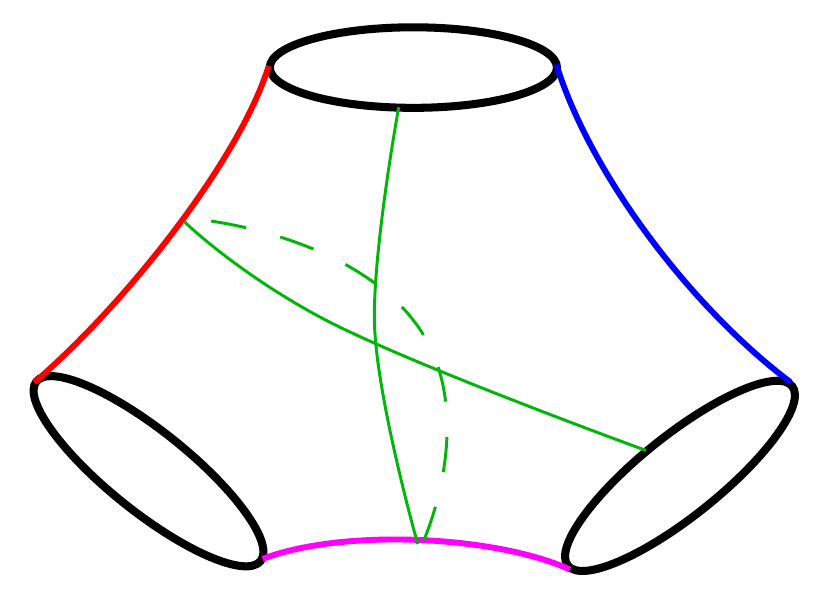}       
    \end{subfigure}
    ~
    \begin{subfigure}{0.25\textwidth}
    \centering
    \includegraphics[{width = \linewidth}]{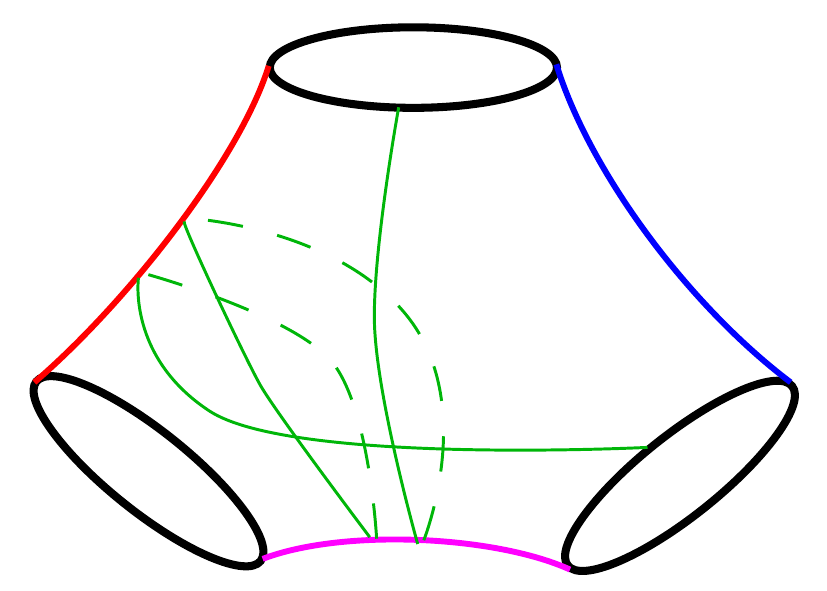}     
    \end{subfigure}
    ~
    \begin{subfigure}{0.25\textwidth}
    \centering
    \includegraphics[{width = \linewidth}]{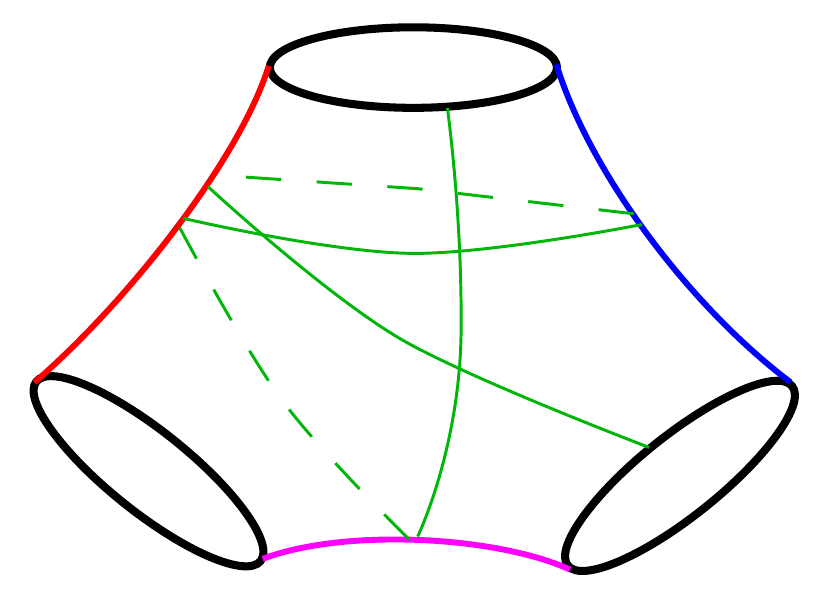}     
    \end{subfigure}
    \newline
    \begin{subfigure}{1\textwidth}
    \centering
    \includegraphics[scale=0.75]{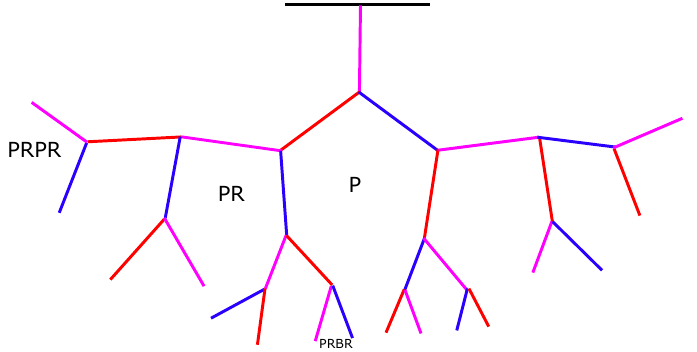}
    \end{subfigure}
    \caption{Consider a pair of pants with a standard orthobasis with labels $\{P,B,R\}$ for pink, blue and red. The top four figures show four orthogeodesic arcs with cutting sequence $P$, $PR$, $PRPR$ and $PRBR$ starting from an interval $I$. The bottom figures show the orthotree parametrizing orthogeodesic starting from $I$.  }
    \label{fig:OrthoTree}
\end{figure}

Orthogeodesics on \( X \in \mathcal{T}(\Sigma_{g,n}) \) starting from the same boundary component \( \delta \) can be encoded by a rooted planar trivalent tree \( T \), known as an \emph{orthotree} \cite[Section~4.1]{doan2021}. The construction of \( T \) depends on a choice of an orthobasis \( \mathcal{H} \) of \( X \) and a choice of a maximal connected open interval $I$ on $\partial X$ disjoint from $\mathcal{H}$. Edges of $T$ are labeled by orthogeodesics in $\mathcal{H}$, and vertices of $T$ are labeled by the hexagons in the decomposition of $X$ by $\mathcal{H}$. We orient all edges away from the root, label the root of $T$ by the initial hexagon that contains $I$ and the edge adjacent to the root of $T$ by the orthogeodesic opposite to $I$ in the initial hexagon. Any finite oriented path starting from the root of $T$ can be thought of as a path on $X$ starting from the initial hexagon crossing a sequence of orthogeodesics in $\mathcal{H}$ transversely. We label the final vertex of the path by the final hexagon in the path on $X$; the left and right edges that point away from this final vertex are labeled by the two remaining left and right orthogeodesics in the final hexagon in the path on $X$. See \cref{fig:OrthoTree} for an example of an orthotree of a pair of pants.   

We observe that each orthogeodesic $\alpha$ starting from $I$ must exit the hexagon containing $I$ through the unique orthogeodesic opposite to $I$. After that, for every hexagon in the decomposition by $\mathcal{H}$ that $\alpha$ enters, there are two possibilities. Either $\alpha$ exits the hexagon by crossing one of the two remaining orthogeodesics or $\alpha$ meets the opposite boundary component perpendicularly. Therefore, each finite oriented path from the root corresponds to an orthogeodesic starting at \( I \subseteq \delta \), crossing orthogeodesics in \( \mathcal{H} \) according to the labels along the path. Thus, orthogeodesics (excluding those in \( \mathcal{H} \)) starting from \( \delta \) correspond to finite paths in \( T \). The following theorem connects the cosh-length spectra of orthogeodesics on pairs of pants and one-holed tori with group orbits.

\begin{thm}[\cite{doan2021}, Theorem~1.2]
\label{thm:GroupAssociatedToOrthoGeodesics}
Let \( P \) be a pair of pants with standard orthobasis \(\mathcal{H}_P = \{\alpha,\beta,\gamma\}\) and three geodesic boundary components $\delta_{\alpha},\delta_{\beta},\delta_{\gamma}$ respectively disjoint from $\alpha,\beta,\gamma$. Let \( T\) be a one-holed torus with orthobasis \(\mathcal{H}_T = \{\alpha,\beta,\gamma\}\). Let $(a,b,c) := (\cosh(\ell(\alpha)),\cosh(\ell(\beta)),\cosh(\ell(\gamma)))$. Then

\begin{enumerate}[label=(\alph*)]
    \item Pair of pants \( P\): The cosh-lengths of orthogeodesics are the positive coordinates of vectors in three orbits \( G_P \cdot \mathbf{u} \), where \( \mathbf{u} \) is one of \((-1,\, c,\, b)^\top\), \((c,\,-1,\, a)^\top\), or \((b,\, a,\,-1)^\top\). Each orbit corresponds to orthogeodesics starting from the boundary components \( \delta_{\alpha},\delta_{\beta},\delta_{\gamma} \) respectively. The group \( G_P \) is generated by:
      \[
    f_{\alpha} = \begin{pmatrix}
        -1 & \frac{2(ab+c)}{a^2-1} & \frac{2(ac+b)}{a^2-1} \\ 
        0 & 1 & 0 \\
        0 & 0 & 1
    \end{pmatrix}, 
    f_{\beta} = \begin{pmatrix}
        1 & 0 & 0 \\ 
        \frac{2(ba+c)}{b^2-1} & -1 & \frac{2(bc+a)}{b^2-1} \\
        0 & 0 & 1
    \end{pmatrix}, 
    f_{\gamma} = \begin{pmatrix}
        1 & 0 & 0 \\ 
        0 & 1 & 0 \\
        \frac{2(ca+b)}{c^2-1} & \frac{2(cb+a)}{c^2-1} & -1
    \end{pmatrix}
   \]
    
    \item One-holed torus \( T \): The cosh-lengths of orthogeodesics are given by the positive coordinates of vectors in the orbit \( G_T \cdot \mathbf{u} \), where \( \mathbf{u} = (b,\, a,\,-1)^\top \). The group \( G_T \) is generated by the matrices:
   \[
    g_{\alpha} = \begin{pmatrix}
        -1 & \frac{b+c}{a-1} & \frac{b+c}{a-1} \\ 
        0 & 0 & 1 \\
        0 & 1 & 0
    \end{pmatrix}, 
    g_{\beta} = \begin{pmatrix}
        0 & 0 & 1 \\ 
        \frac{c+a}{b-1} & -1 & \frac{c+a}{b-1} \\
        1 & 0 & 0
    \end{pmatrix}, 
    g_{\gamma} = \begin{pmatrix}
        0 & 1 & 0 \\ 
        1 & 0 & 0 \\
        \frac{a+b}{c-1} & \frac{a+b}{c-1} & -1
    \end{pmatrix}.
    \]
\end{enumerate}
\end{thm}

\begin{remark}
    Two orthogeodesics are called \emph{neighbors} with respect to an orthobasis \( \mathcal{H} \) if they originate from the same boundary component and, together with an element of \( \mathcal{H} \), form a triple of alternating sides of an immersed hexagon on the surface. Moreover, each vector with positive coordinates in the orbits of \cref{thm:GroupAssociatedToOrthoGeodesics} corresponds to the cosh-lengths of three pairwise neighboring orthogeodesics.
\end{remark}

We recall some result about short orthogeodesics on $X$. The \emph{orthosystole}, \(\operatorname{osys}(X)\), of a hyperbolic surface \(X\) is the length of its shortest orthogeodesic. The \emph{orthokissing number}, $\operatorname{okiss}(X)$, is the number
of orthogeodesics of length \(\operatorname{osys}(X)\). Bavard \cite{bavard2005anneaux} established a sharp upper bound for the orthosystole of a hyperbolic surface with totally geodesic boundary, restated equivalently by Basmajian and Fanoni \cite[Theorem~1.1]{basmajian2024orthosystoles}. The following is an equivalent formulation:

\begin{thm}[Bavard]\label{thm:Bavardupperbound}
Let \( X \) be a hyperbolic surface of signature \( (g, n) \) with total boundary length \( L \). Then
\[
\cosh\bigl( \osys(X) \bigr) \leq \frac{1}{\cosh\left( \dfrac{L}{12g - 12 + 6n} \right) - 1} + 1.
\]
 Equality holds if and only if \( \operatorname{okiss}(X) = 6g - 6 + 3n \). Moreover, the bound is attained for every \( L \).
\end{thm}

\begin{remark}
Bavard’s original work is more general and presented differently. The theorem above corresponds to parts (1), (3), and (4) of \cite[Théorème~1]{bavard2005anneaux}, with \( c = \partial X \) and \( r = \dfrac{\osys(X)}{2} \). In this context, the equality condition in part (3) becomes \( \operatorname{okiss}(X) = 6g - 6 + 3n \). For details on the correspondence between Bavard's formulation and the theorem above, see \cite[Appendix]{basmajian2024orthosystoles}.
\end{remark}

\subsection{Arithmetic Fuchsian groups coming from quadratic forms}

The material presented in this section is discussed in \cite[Section 10.1]{maclachlan2003book}, and will be used exclusively in \cref{sec:Ortho-integral surfaces and quadratic forms}. Let $V$ be a 3-dimensional vector space over $\mathbb{R}$ equipped with a quadratic form $q$ of signature $(2,1)$. Therefore, there exists a basis $\{\mathbf{e}_1, \mathbf{e}_2, \mathbf{e}_3\}$ for $V$ such that
\[
q(\mathbf{x}) = x_1^2 + x_2^2 - x_3^2,
\]
where $\mathbf{x} = x_1 \mathbf{e}_1 + x_2 \mathbf{e}_2 + x_3 \mathbf{e}_3 \in V$.

The hyperbolic plane $\mathbb{H}^2$ can be modeled by the \emph{hyperboloid model}, identifying $\mathbb{H}^2$ with the upper sheet
\[
C^+ = \{ \mathbf{x} \in V  \mid q(\mathbf{x}) = - 1 \text{ and } x_3 >0 \}.
\]
Under this identification, the isometry group of \( \mathbb{H}^2 \) is given by
\[
\mathrm{O}^+(V,q) = \{ g  \in \mathrm{O}(V,q) \mid g(C^+) = C^+ \},
\]
which is the subgroup of the orthogonal group \( \mathrm{O}(V, q) \) that preserves \( C^+ \).

Arithmetically defined subgroups of $\Isom(\mathbb{H}^2)$ arise from quadratic forms defined over number fields. Let $k \subset \mathbb{R}$ be a number field, and let $V$ be a 3-dimensional vector space over $k$ equipped with a quadratic form $q$ defined over $k$ and of signature $(2,1)$ over $\mathbb{R}$. By extending scalars, we obtain $V \otimes_k \mathbb{R}$, and $q$ extends naturally to this space. The orthogonal group $\mathrm{O}(V, q)$ over $k$ embeds into the real orthogonal group:
\[
\mathrm{O}(V, q)(k) \hookrightarrow \mathrm{O}(V \otimes_k \mathbb{R}, q).
\]
Let \( L \) be a \emph{lattice} in \( V \), defined as a free \( \mathcal{O}_k \)-module of rank three, where \( \mathcal{O}_k \) is the ring of integers in \( k \). The \emph{special orthogonal group} associated to \( L \) is
\[
\mathrm{SO}(L) = \left\{ g \in \mathrm{SO}^+(V, q)(k) \;\middle|\; g(L) = L \right\}.
\]
Here, \( \mathrm{SO}^+(V, q)(k) \) denotes the subgroup of \( \mathrm{O}^+(V, q)(k) \) consisting of elements with determinant 1 that preserve the orientation of \( C^+ \).

We are interested in computing the invariant quaternion algebra of $\SO(L)$ from the data of the quadratic space $(V,q)$ over $k$. The connection between $\SO(L)$ and discrete subgroups of $\PSL_2(\mathbb{R})$ and $\PSL_2(\mathbb{C})$  can be established using Clifford algebras. The Clifford algebra $C(V)$ of a non-degenerate quadratic space $(V,q)$ over a field $k$ is an associative algebra with $1$ containing $V$ such that $x^2 = q(x)1$ for all $x \in V$. When $V$ is 4-dimensional over $k$, the algebra $C(V)$ is a $\mathbb{Z}/2\mathbb{Z}$-graded algebra of dimension 16 over $k$. The even sub-algebra $C_0(V)$ has dimension 8 over $k$ and  is isomorphic to a quaternion algebra over \( k \):
\begin{equation}
\label{eq:QuaternionAlgebraAssociatedToQuadraticSpace}
    C_0(V) \cong \left( \frac{ -q(\mathbf{e}_1) q(\mathbf{e}_2),\ -q(\mathbf{e}_1) q(\mathbf{e}_3) }{ k } \right),
\end{equation}  
where $\{\mathbf{e}_1, \mathbf{e}_2, \mathbf{e}_3\}$ is an orthogonal basis for $V$ such that $q({\bf e}_1) < 0$ . In \cref{sec:Ortho-integral surfaces and quadratic forms}, we establish the connection between OI pair of pants and OI one-holed torus and a certain quadratic form that depends on a choice of an orthobasis. This connection allows a straightforward verification of arithmeticity of the double of these surfaces along their boundary components.

\subsection{Hyperbolic identities}

From \cite[Theorem~2.4.1]{buser2010geometry}, any hyperbolic convex right-angled geodesic hexagon with consecutive sides \(\alpha, \beta', \gamma, \alpha', \beta, \gamma'\) satisfies:
\begin{equation}
    \label{eq:RightAngledHexagonIdentities}
    \cosh{\ell(\gamma')} = \frac{\cosh{\ell(\alpha)}\cosh{\ell(\beta)} + \cosh{\ell(\gamma)}}{\sinh{\ell(\alpha)}\sinh{\ell(\beta)}}.
\end{equation}
\noindent
 
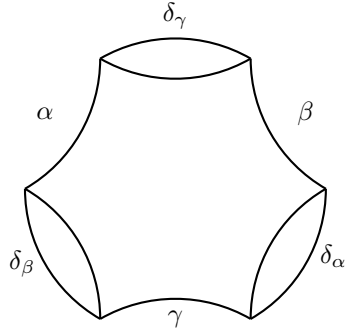
\begin{figure}[h!]
    \centering
    \begin{tikzpicture}[scale = 2]
        % Draw the complete outline of the pair of pants.
        \draw[thick] (1,0) arc (240:180:1);
        \draw[thick] (1/2,{sqrt(3)/2}) arc (300:240:1);
        \draw[thick] (-1/2,{sqrt(3)/2}) arc (0:-60:1);
        \draw[thick] (-1,0) arc (60:0:1);
        \draw[thick] (-1/2,{-sqrt(3)/2}) arc (120:60:1);
        \draw[thick] (1/2,{-sqrt(3)/2}) arc (180:120:1);
        %%%
        % Draw and label the three boundary arcs.
        \draw[thick] (1/2,{sqrt(3)/2}) arc (60:120:1) node[pos=0.5, anchor=south] {$\delta_{\gamma}$};
        \draw[thick] (-1,0) arc (180:240:1) node[pos=0.5, anchor=east] {$\delta_{\beta}$};
        \draw[thick] (1/2,{-sqrt(3)/2}) arc (300:360:1) node[pos=0.7, anchor=north, xshift=2mm] {$\delta_{\alpha}$};
        %%%
        % Label the orthogeodesics.
        \draw ({sqrt(3)/2},1/2) node {$\beta$};
        \draw ({-sqrt(3)/2},1/2) node {$\alpha$};
        \draw (0,{-sqrt(3)/2}) node {$\gamma$};
    \end{tikzpicture}  
    \caption{A pair of pants \(P\) with the standard orthobasis \(\mathcal{H} = \{\alpha,\beta,\gamma\}\) and three geodesic boundary components $\delta_{\alpha},\delta_{\beta},$ and $\delta_{\gamma}$.}
    \label{fig:Pants}
\end{figure}

Throughout this paper, \( P \) is always a pair of pants with standard orthobasis \(\mathcal{H}_P = \{\alpha,\beta,\gamma\}\) and three geodesic boundary components $\delta_{\alpha},\delta_{\beta},\delta_{\gamma}$ respectively disjoint from $\alpha,\beta,\gamma$ (see \cref{fig:Pants}). Denote $(a,b,c) := (\cosh(\ell(\alpha)),\cosh(\ell(\beta)),\cosh(\ell(\gamma)))$.

\begin{lem}
\label{lem:TraceOfBoundaries}
 We have
\begin{equation}
\label{eq:FormulasForTraceSquares}
\tr_{\SL}^2(\delta_\alpha) = \frac{4(bc+a)^2}{(b^2-1)(c^2-1)}, \tr_{\SL}^2(\delta_\beta) = \frac{4(ca+b)^2}{(c^2-1)(a^2-1)}, \tr_{\SL}^2(\delta_\gamma) = \frac{4(ab+c)^2}{(a^2-1)(b^2-1)}.
\end{equation}

\end{lem}

\begin{proof}
It suffices to check the identity for $\tr_{\SL}^2(\delta_\alpha)$ since the other two can be checked in a similar manner. Applying \cref{eq:RightAngledHexagonIdentities} to the hexagon with alternating edges $\alpha$, $\beta$ and $\gamma$, we get
\[
\cosh\left(\frac{\ell(\delta_\alpha)}{2}\right) = \frac{(bc+a)}{\sqrt{(b^2-1)(c^2-1)}}.
\]
The lemma follows since $\tr^2_{\SL}(\delta_\alpha) = 4\cosh^2\left(\frac{\ell(\delta_\alpha)}{2}\right)$.

    \end{proof}

\begin{lem}
\label{lem:CoshLengthOrthoGeodesicAroundBoundary}
 Let $\eta_k$ be the oriented orthogeodesic with cutting sequence of length $2k+1$ of the form $(\gamma,\beta, \dots,\gamma)$ where $k\geq 0$. The cosh-length of the orthogeodesic $\eta_{k}$ can be computed recursively by 
\begin{equation}
\label{eq:RecursionCoshLengthOrthoGeodesicAroundBoundary}
    \cosh(\ell(\eta_{k+3})) = (\tau-1)\cosh(\ell(\eta_{k+2})) - (\tau-1)\cosh(\ell(\eta_{k+1})) + \cosh(\ell(\eta_{k}))
\end{equation}
with initial conditions $\cosh(\ell(\eta_{0})) = s-1$, $\cosh(\ell(\eta_{1})) = \tau s-1$ and $\cosh(\ell(\eta_{2})) = \tau^2s-2\tau s +s-1$ where
\begin{equation*}
    \tau = \tr_{\SL}^2(\delta_\alpha) = \frac{4(bc+a)^2}{(b^2-1)(c^2-1)}\quad \text{and}\quad s = \frac{2 ( a^2 + b^2 + c^2 + 2abc -1)}{c^2-1}.
\end{equation*}
\end{lem}

\begin{proof}
    We note that $f_\beta$ and $f_\gamma$ (see \cref{thm:GroupAssociatedToOrthoGeodesics}) fixes the vector $(-1,c,b)^\top$ and that $$\tr_{\SO}(f_\gamma f_\beta) = \frac{4(bc+a)^2}{(b^2-1)(c^2-1)}-1 =\tr_{\SL}^2(\delta_\alpha)-1.$$ For brevity, we write $\tau = \tr_{\SL}^2(\delta_\alpha)$ and $\textbf{v}= (b,a,-1)^\top$. The characteristic polynomial of $f_\gamma f_\beta$ has the form 
    \[\chi(\lambda) = 
    (\lambda-1)(\lambda^2 + (2-\tau)\lambda+1) =  \lambda^3 + (1-\tau)\lambda^2 +( \tau-1)\lambda - 1.
    \]
    By the Cayley–Hamilton theorem, $f_\gamma f_\beta$ satisfies its own characteristic polynomial, we have
    \[
    (f_\gamma f_\beta)^{k+3} = (\tau - 1)(f_\gamma f_\beta)^{k+2} -  (\tau - 1)(f_\gamma f_\beta)^{k+1} +  (f_\gamma f_\beta)^{k}.
    \]
    Since $\cosh(\ell(\eta_k))$ is the last coordinate of $(f_\gamma f_\beta)^{k}f_\gamma \textbf{v}$, the recursive formula for $\cosh(\ell(\eta_k))$ follows from the recursive formula for $ (f_\gamma f_\beta)^{k}$.
\end{proof}

\section{The finiteness of OI surfaces}\label{sec:Finiteness}

We now introduce a broader class of surfaces, which we call 
$d$-AOI surfaces, and prove a finiteness result for them that implies the finiteness of OI surfaces in \cref{thm:FinitenessOfOISurfaces}.
\begin{defi}
    \label{def:OIandAlmostOI}
  
  A hyperbolic surface \( X \) is  \emph{\( d \)-almost-ortho-integral (\( d \)-AOI)} for some \( d \in \mathbb{Z}^{+} \) if the cosh of the length of every orthogeodesic of \( X \) is a number in \( \frac{1}{d}\mathbb{Z} \).

\end{defi}

In this section, we provide the proof of the following theorem.

\begin{thm}\label{thm:FinitenessOfdAOISurfaces}
    For all $d\in \mathbb{Z}^{+}$, the set of $d$-AOI surfaces of a fixed genus $g$ and $n>0$ boundary components is finite. 
\end{thm}

We begin our discussion with the following definition.

\begin{defi}
  Let \( X \) be a hyperbolic surface with nonempty totally geodesic boundary. Fix a boundary component \(\delta\) of \(X\). We say that \(X\) is \emph{partially \(d\)-AOI} with respect to \(\delta\) if the cosh of the length of every orthogeodesic with both endpoints on \(\delta\) lies in \(\frac{1}{d}\mathbb{Z}\). 

\end{defi}

\begin{lem}\label{lem:partiallyOI}
Let $P$ be a hyperbolic pair of pants with a standard orthobasis $\mathcal{H} = \{\alpha,\beta,\gamma\}$. Denote $(a,b,c) := (\cosh\ell(\alpha),\cosh\ell(\beta),\cosh\ell(\gamma))$. If $P$ is partially $d$-AOI with respect to the boundary $\delta_{\gamma}$, then 
\begin{equation}
\tr_{\SL}^2(\delta_\alpha) = \frac{4(bc+a)^2}{(b^2-1)(c^2-1)}, \quad \text{and}\quad \tr_{\SL}^2(\delta_\beta) = \frac{4(ca+b)^2}{(c^2-1)(a^2-1)}
\end{equation}
are all positive integers.
\end{lem}

\begin{proof}
    It suffices to prove that $\tr_{\SL}^2(\delta_\alpha) \in \mathbb{N}$. The other case followed by a similar argument. Let $\eta_k$ be the oriented orthogeodesic with cutting sequence of length $2k+1$ of the form $(\gamma,\beta, \dots,\gamma)$ where $k\geq 0$. By \cref{lem:CoshLengthOrthoGeodesicAroundBoundary}, we have $\cosh(\ell(\eta_{0})) = s-1$, $\cosh(\ell(\eta_{1})) = \tau s-1$ and $\cosh(\ell(\eta_{2})) = \tau^2s-2\tau s +s-1$ where
    \begin{equation}
    \label{eq:InitialConditions}
    \tau = \tr_{\SL}^2(\delta_\alpha) = \frac{4(bc+a)^2}{(b^2-1)(c^2-1)}\quad \text{and}\quad s = \frac{2 ( a^2 + b^2 + c^2 + 2abc -1)}{c^2-1}.
    \end{equation}
    \noindent
    Since $P$ is partially $d$-AOI with respect to the boundary $\delta_{\gamma}$, $\cosh(\ell(\eta_{k})) \in \frac{1}{d}\mathbb{Z}$ for all $k \in \mathbb{N}$. The initial conditions for $\cosh(\ell(\eta_{k}))$ implies that $\{s,\tau s, \tau^2 s\} \subset \frac{1}{d}\mathbb{Z}$. From the initial conditions and the recursive formula for $\cosh(\ell(\eta_k))$ in \cref{eq:RecursionCoshLengthOrthoGeodesicAroundBoundary}, the formula for $\cosh(\ell(\eta_k))$ is a $\mathbb{Z}$-linear combination of $\{\tau^k s,\tau^{k-1}s\dots, s,1\}$ where the coefficient of $\tau^k s$ is one. Arguing by induction, we get that $\tau^k s \in \frac{1}{d}\mathbb{Z}$ for all integers $k \geq 0$. We see that $\tau \in \mathbb{Q}$ since $s,\tau s \in \frac{1}{d}\mathbb{Z}$. Thus, $\tau^k s \in \frac{1}{d}\mathbb{Z}$ for all integers $k \geq 0$ if and only if $\tau = \tr_{\SL}^2(\delta_\alpha) \in \mathbb{N}$. 
\end{proof}

The following proposition provides the key bounds needed to establish the finiteness of the set of \(d\)--AOI surfaces. The crucial observation is that for any two (not necessary distinct) closed geodesics on a hyperbolic surface \(X\), there exists an immersed pair of pants bounded by these two geodesics. This pair of pants can be constructed explicitly as follows. Let $\Gamma\leq \PSL_2(\mathbb{R})$ be the image of $\pi_1(X)$ under the holonomy representation. For any closed geodesics $\mu_{1}$ and $\mu_2$ on $X$, we choose two disjoint lifts $\widetilde{\mu}_1$ of $\mu_1$ and $\widetilde{\mu}_2$ of $\mu_2$ in the universal cover $\widetilde{X} \subset \mathbb{H}^2$. The isometries $g_{\mu_1} \in \Gamma$ stabilizing $\widetilde{\mu}_1$ and $g_{\mu_2} \in \Gamma$ stabilizing $\widetilde{\mu}_2$ generate a subgroup $\Delta$ in $\Gamma$. Let $C(\Lambda_\Delta)$ be the convex hull of the limit set of $\Delta$. Then the quotient $C(\Lambda_\Delta)/\Delta$ is isometric to a pair of pants $P$. The natural inclusion of $C(\Lambda_\Delta)$ in $\widetilde{X} \subset \mathbb{H}^2$ induces an isometric immersion of the pair of pants $P$ into $X$. 

\begin{remark}
Immersed pairs of pants were used in \cite{basmajian2020prime, parlier2020geodesic} to extend Basmajian's identity beyond hyperbolic surfaces with geodesic boundary.
\end{remark}
\begin{prop}\label{prop:IntegralityOfdAOISurfaces}
    Let $X\in \mathcal{M}_{g,n}$ be a $d$-AOI surface. Let $L$ be the total boundary length of $X$. Then  
    \begin{itemize}
    \item[(i)] $\tr_{\SL}^2(\mu)\in \mathbb{Z}$ for any closed geodesic $\mu$ on $X$.
    \item[(ii)] $L\leq (12g-12+6n)\operatorname{arcosh}(d+1).$
    \item[(iii)] $\operatorname{sys}(X) \geq 2\operatorname{arcosh}\left(\frac{\sqrt{5}}{2}\right)$.
    \item[(iv)] $\ell(\delta) \geq 2\operatorname{arcosh}\left(\frac{\sqrt{5}}{2}\right)$, for any boundary component $\delta$ of $X$.
    \end{itemize}
\end{prop}

\begin{proof}
$(i)$ Let $P$ be an immersed pair of pants on $X$ with three geodesic boundary components $\delta_{\alpha},\delta_{\beta},\delta_{\gamma}, $ where $\delta_{\alpha}=\mu$ and $\delta_{\gamma}  \in \partial X$. The pair of pants $P$ is constructed as described in the preceding paragraph. Since $X$ is a $d$--AOI surface, $P$ is partially $d$--AOI with respect to $\delta_{\gamma}$. By \cref{lem:partiallyOI}, $\tr_{\SL}^2(\mu)\in \mathbb{Z}$.\\
$(ii)$ Since $X$ is $d$--AOI, $\cosh(\osys(X))$ is greater than or equal to $1+\frac{1}{d}$. Applying \cref{thm:Bavardupperbound}, we have
\[\cosh\left( \dfrac{L}{12g - 12 + 6n} \right) \le \frac{1}{\cosh\bigl( \osys(P) \bigr)-1}+1 \le d+1,
\]
which yields the upper bound on $L$. We observe: $(i)$ implies that $\tr_{\SL}(\mu)^2 \geq 5$ for any closed geodesic on $X$. Therefore, $(iii)$ and $(iv)$ follows from $(i)$ and the fact that $\cosh(\ell(\mu)/2) = \frac{|\tr_{\SL}(\mu)|}{2}$. 
\end{proof}

\begin{remark}\label{rem:IntegralityLargerClass}
    Properties $(i)$, $(iii)$, and $(iv)$ in \cref{prop:IntegralityOfdAOISurfaces} still hold if \( X \) is a partially \( d \)-AOI surface with respect to some fixed boundary component \( \delta \).
\end{remark}

\begin{proof}[Proof of \cref{thm:FinitenessOfdAOISurfaces}]
 \cref{prop:IntegralityOfdAOISurfaces} provides uniform upper and lower bounds on the length of each boundary component of $X$ and establishes a lower bound on the systole of \( X \). By Mumford's compactness criterion \cref{thm:MumfordCompactnesCriterion}, the set of $d$--AOI  surfaces \( X \in \mathcal{M}_{g,n}\) is contained within a compact subset of the moduli space \( \mathcal{M}_{g,n} \).

Since the cosh-lengths of orthogeodesics of $d$--AOI surfaces can only take value in a fixed discrete subset of $\mathbb{R}$, the set of $d$--AOI surfaces is discrete in the moduli space $\mathcal{M}_{g,n}$. Together with the compactness, we conclude that there are only finitely many $d$--AOI surfaces for a fixed $d$ and fixed topological complexity. 

\end{proof}

\begin{remark}\label{rem:DiscretenessPAOI}
The discreteness in the moduli space \(\mathcal{M}_{g,n}\) still holds for the set of partially \(d\)--AOI surfaces since $\tr_{\SL}^2(\mu) \in \mathbb{Z}$ for every closed geodesic $\mu$ on these surfaces. However, finiteness may fail for this broader class, as we do not have a uniform lower bound on the orthosystole, which in turn prevents us from deriving a uniform upper bound on the total boundary length.
\end{remark}

\section{Classification of OI pairs of pants and OI one-holed tori}\label{sec:Classification of OI pants and OI one-holed tori}

\subsection{Classification of OI pairs of pants}

The following result provides another way to prove the finiteness of ortho-integral pairs of pants without using the Mumford's compactness criterion and offers a complete classification of them.

\begin{thm}
\label{thm:classificationOIpants}
Let $P(a_1,a_2,a_3)$ be a pair of pants with the standard orthobasis $\{\alpha, \beta, \gamma\}$. Denote $(a_1,a_2,a_3) := (\cosh(\ell(\alpha)),\cosh(\ell(\beta)),\cosh(\ell(\gamma)))$. Assume that $a_1 \le a_2 \le a_3$, then the following are equivalent:
\begin{itemize}
    
    \item[(i)] $P(a_1,a_2,a_3)$ is ortho-integral.
   
    \item[(ii)]  $(a_1,a_2,a_3)\in \{(2,2,2),(2,2,5),(2,2,17),(3,3,3),(3,3,7),(3,3,19),(5,5,11)\}.$

     \item[(iii)] $\left\{a_i,\dfrac{2(a_1^2+a_2^2+a_3^2+2a_1a_2a_3-1)}{a_i^2-1},\dfrac{4(a_ia_j+a_k)^2}{(a_i^2-1)(a_j^2-1)}\right\}\subset \mathbb{N},$ for $\{i,j,k\} = \{1,2,3\}$.

\end{itemize}
\end{thm}
\begin{proof}
$(ii) \Rightarrow (i):$ The proof is given in \cite[Theorem 1.1]{doan2021}.\\
$(i) \Rightarrow (iii):$ Directly from \cref{lem:partiallyOI} and its proof when $d=1$.\\
$(iii) \Rightarrow (ii):$ If $a_2 \ge 3$ and $a_3>19$, then $$4+\frac{4(a_1^2+a_3^2+2a_1a_2a_3)}{a_2^2-1} \leq 4+\frac{4(a_2^2+a_3^2+2a_2^2a_3)}{a_2^2-1} = 4+4\left(2a_3+1+\frac{(a_3+1)^2}{a_2^2-1}\right)$$
    $$\le 4+4\left(2a_3+1+\frac{(a_3+1)^2}{8}\right) = \frac{a_3^2}{2}+9a_3+\frac{17}{2}<a_3^2-1.$$
 It implies that
\[
\frac{4(a_2 a_3 + a_1)^2}{(a_2^2 - 1)(a_3^2 - 1)} = 4 + \frac{4 + \frac{4(a_1^2 + a_3^2 + 2a_1 a_2 a_3)}{a_2^2 - 1}}{a_3^2 - 1} \notin \mathbb{N},
\]
which contradicts condition~$(iii)$. Therefore, either \( a_2 = 2 \) or \( a_3 \leq 19 \). 

If $a_2 = 2$, then $2=a_1=a_2 \le a_3$, and 
$$\frac{2(a_1^2+a_2^2+a_3^2+2a_1a_2a_3-1)}{a_3^2-1} = \frac{2(4+4+a_3^2+8a_3-1)}{a_3^2-1} = 2+\frac{16}{a_3-1}.$$ By $(iii)$, $(a_3-1)|16$, thus $a_3 \in \{2,3,5,9,17\}$. One can check that condition~$(iii)$ satisfies only if $a_3 \in \{2,5,17\}$.

If $a_3 \le 19$, then $2 \le a_1 \le a_2 \le a_3 \le 19$. One can verify using a short computer program that the list provided in $(ii)$ contains all possible solutions satisfying condition~$(iii)$.    
\end{proof}

\subsection{Classification of OI one-holed tori}

In this subsection, we will classify the set of OI one-holed tori. Before stating the result, let us make the following definition.

\begin{defi}\label{def:Minimalorthobasis}
An orthobasis \(\{\alpha,\beta,\gamma\}\) of a one-holed torus \(T\) is called \emph{minimal} if it contains one of the shortest orthogeodesics on \(T\) (an orthosystole geodesic) and if 
\[
\max\{\ell(\alpha),\ell(\beta),\ell(\gamma)\}
\]
is the smallest among all orthobases that contain an orthosystole geodesic.
\end{defi}

\begin{remark}
   Since the ortholength spectrum of a surface is discrete in $\mathbb{R}$, a minimal orthobasis necessarily exists, and there can be only finitely many minimal orthobases. 
\end{remark}

Moreover, as we will prove in \cref{lem:propertyofMinimal}, the cosh-length triple associated with any minimal orthobasis on a given one-holed torus is unique up to permutation. Consequently, the one-holed tori can be uniquely parametrized up to isometry by their minimal orthobasis cosh-length triples. The following theorem presents the classification up to isometry of OI one-holed tori.

\begin{thm}\label{thm:classificationOItorus}
Let $(a,b,c),$ where $ 1 < a\le b \le c$, be the triple of cosh-lengths of a minimal orthobasis of a one-holed torus $T$, then the following are equivalent:
\begin{itemize}

    \item[(i)] $T$ is ortho-integral.
    \item[(ii)]  $(a,b,c)\in \{(2, 2, 2),(2, 2, 3),(2, 2, 5),(2, 3, 6),(2, 4, 4),(2, 4, 7),(2, 5, 8),(2, 7, 10),(2, 13, 16),\\(3, 3, 3),(3, 3, 7),(3, 5, 5),(3, 5, 9),(3, 9, 13),(3, 17, 21),(4, 4, 5),(4, 5, 10),(4, 6, 6),(4, 11, 16),\\(5, 5, 11),(5, 7, 13),(5, 13, 19),(6, 6, 9),(6, 36, 64),(7, 9, 9),(7, 17, 25),(8, 13, 22),(9, 11, 21),\\(10, 10, 17),(10, 12, 12),(11, 49, 61),(13, 29, 43),(17, 19, 37),(19, 21, 21)\}.$

     \item[(iii)] $2\le a \le 20$, $c-a-1 \le b\le c \le 107$, and $$\left\{a,b,c,\frac{(a+b+c-1)^2}{(a-1)(b-1)},\frac{(a+b+c-1)^2}{(b-1)(c-1)},\frac{(a+b+c-1)^2}{(c-1)(a-1)}\right\}\subset \mathbb{N}.$$

\end{itemize}
\end{thm}

For $\alpha$ be a simple orthogeodesic on a one-holed torus $T$, the \emph{accompanying simple closed geodesic} $\lambda_{\alpha}$ is the unique simple closed geodesic on $T$ disjoint from $\alpha$.

\begin{lem}\label{lem:IdentitiesOnOneHoledTorus}
    Let $(a,b,c) := (\cosh\ell(\alpha),\cosh\ell(\beta),\cosh\ell(\gamma))$ denote the triple of the cosh-lengths of an orthobasis $\{\alpha,\beta,\gamma\}$ for a one-holed torus \( T \), and let $ (\lambda_\alpha, \lambda_\beta, \lambda_\gamma) $ be the accompanying triple of simple closed geodesics. Then
\begin{equation}\label{eq:BoundaryOneHoledTorus}
    |\tr_{\SL}(\partial T)| =  2\cosh(L/2)= \frac{2(a+b+c-1)^2}{(a-1)(b-1)(c-1)}+2, \text{ and }
\end{equation}  
\begin{equation*}
    \tr_{\SL}^2(\lambda_\alpha)=\frac{(a+b+c-1)^2}{(b-1)(c-1)}, \ \tr_{\SL}^2(\lambda_\beta)=\frac{(a+b+c-1)^2}{(c-1)(a-1)}, \
    \tr_{\SL}^2(\lambda_\gamma)=\frac{(a+b+c-1)^2}{(a-1)(b-1)}.
\end{equation*}
where $L = \ell(\partial T)$. Furthermore, if $T$ is an OI one-holed torus then $\tr_{\SL}(\partial T) \in \mathbb{Z}$.  
\end{lem}
\begin{proof} Let $\eta$ and $\xi$ be orthogeodesics on $T$ with the corresponding cutting sequences $\{\gamma\}$ and $\{\gamma,\beta\}$ respectively. For brevity, we write $p = \cosh(\ell(\eta))$ and $q=\cosh(\ell(\xi))$. By \cref{thm:GroupAssociatedToOrthoGeodesics}, $g_{\gamma}((b,a,-1)^\top)=(a,b,p)$ and $g_{\beta}((a,b,p)^\top)=(p,q,a)$, thus
    $$ p = 1+\frac{a+b}{c-1}\cdot (a+b),$$
    and $$ q = -b+\frac{a+c}{b-1}\cdot (p+a) = -b+\frac{a+c}{b-1}\cdot\left(1+\frac{(a+b)^2}{c-1}+a\right).$$
We observe that \(T\) contains an immersed right-angled hexagon with alternating sides \((\alpha, \xi, \alpha)\). Lifting this hexagon to the universal cover of \(T\) yields a right-angled hexagon with alternating sides whose cosh-lengths are \((a, q, a)\), and the side opposite the lift of \(\xi\) has length \(L/2\). By the identity \cref{eq:RightAngledHexagonIdentities} for the right-angle hexagon, we have that $$ 2\cosh(L/2) = 2\frac{a^2+q}{a^2-1}
= 2\frac{a^2-b+\frac{a+c}{b-1}\left(1+\frac{(a+b)^2}{c-1}+a\right)}{a^2-1}
= \frac{2(a+b+c-1)^2}{(a-1)(b-1)(c-1)}+2.$$

Now, we prove the identities for the accompanying geodesics. Cutting $T$ along $\alpha$ and the accompanying geodesic $\lambda_\alpha$ yields two isometric half-pants. Each half-pants can be decomposed into two copies of a right-angled pentagon with side length $L/4$ adjacent to $\ell(\alpha)/2$ and with side length $\ell(\lambda_\alpha)/2$ disjoint from both $L/4$,  and $\ell(\alpha)/2$. Using standard hyperbolic trigonometry \cite[Theorem 2.3.4]{buser2010geometry}, we have
 $$\cosh \frac{\ell(\lambda_\alpha)}{2} = \sinh \frac{\ell(\alpha)}{2}\sinh\frac{L}{4}.$$
 Combining with the identity for $|\tr_{\SL}(\partial T)|$, we have $$\tr_{\SL}^2(\lambda_\alpha) = (\cosh \ell(\alpha)-1)(\cosh(L/2)-1)=(a-1)\cdot \frac{(a+b+c-1)^2}{(a-1)(b-1)(c-1)}= \frac{(a+b+c-1)^2}{(b-1)(c-1)}.$$
The other cases follow analogously. 

Now, we assume that $T$ is an OI one-holed torus. \cref{prop:IntegralityOfdAOISurfaces} and the previous identities imply that 
\begin{equation*}
    \tr_{\SL}^2(\lambda_\alpha)=\frac{(a+b+c-1)^2}{(b-1)(c-1)}, \ \tr_{\SL}^2(\lambda_\beta)=\frac{(a+b+c-1)^2}{(c-1)(a-1)}, \
    \tr_{\SL}^2(\lambda_\gamma)=\frac{(a+b+c-1)^2}{(a-1)(b-1)}
\end{equation*}
are all integers. This implies that 
    \begin{equation}\label{equ:squareproduct}
        (\tr_{\SL}(\lambda_\alpha)\tr_{\SL}(\lambda_\beta)\tr_{\SL}(\lambda_\gamma))^2=\left(\frac{(a+b+c-1)^3}{(a-1)(b-1)(c-1)}\right)^2\in \mathbb{N}.
    \end{equation} 
    On the other hand since $a,b,c \in \mathbb{N}$, we have
     \begin{equation}\label{equ:product}
         |\tr_{\SL}(\lambda_\alpha)\tr_{\SL}(\lambda_\beta)\tr_{\SL}(\lambda_\gamma)|=\frac{(a+b+c-1)^3}{(a-1)(b-1)(c-1)}\in \mathbb{Q}.
     \end{equation} \cref{equ:squareproduct} and \cref{equ:product} imply that $\tr_{\SL}(\lambda_\alpha)\tr_{\SL}(\lambda_\beta)\tr_{\SL}(\lambda_\gamma)\in \mathbb{Z}$. The Markoff-type equation for one-holed tori
$$\tr_{\SL}^2(\lambda_\alpha)+\tr_{\SL}^2(\lambda_\beta)+\tr_{\SL}^2(\lambda_\gamma)-\tr_{\SL}(\lambda_\alpha)\tr_{\SL}(\lambda_\beta)\tr_{\SL}(\lambda_\gamma) = 2-|\tr_{\SL}(\partial T)|$$
implies that $|\tr_{\SL}(\partial T)|\in \mathbb{N}$.
\end{proof}

\begin{lem}\label{lem:propertyofMinimal}
    Let $\{\alpha,\beta,\gamma\}$ be a minimal orthobasis of a one-holed torus $T$ with $$(a,b,c) := (\cosh\ell(\alpha),\cosh\ell(\beta),\cosh\ell(\gamma)),$$ where $ 1 < a\le b \le c$. Then
    \begin{enumerate}
        \item[(i)] $a =\cosh(\operatorname{osys}(T)) \ge c-b-1,$

        \item[(ii)]If $(a,b',c)$, where $a\le b'\le c$, is the cosh-length triple of another minimal orthobasis of $T$, then $b'= b$.
    \end{enumerate}
    
\end{lem}
\begin{proof} $(i)$ The fact that $a =\cosh(\operatorname{osys}(T))$ follows immediately from \cref{def:Minimalorthobasis}. Let \(\gamma'\) be the simple orthogeodesic with cutting sequence \(\{\gamma\}\). By applying \cref{thm:GroupAssociatedToOrthoGeodesics}, we obtain
\[
\cosh\ell(\gamma') = \frac{(a+b)^2}{\,c-1}+1,
\]
which shows that the cosh-length triple corresponding to the orthobasis \(\{\alpha,\beta,\gamma'\}\) is
\[
\left(a,\, b,\, \frac{(a+b)^2}{\,c-1}+1\right).
\]
Since the original orthobasis \(\{\alpha,\beta,\gamma\}\) is minimal, we must have
\[
c \le \max\left\{b, \frac{(a+b)^2}{\,c-1}+1\right\}.
\]
\begin{itemize}
    \item If $c \le \max\left\{b, \frac{(a+b)^2}{\,c-1}+1\right\}=b$, then $b=c$, and $$b \ge \frac{(a+b)^2}{\,c-1}+1 =\frac{a+c}{c-1}(a+b)+1>a+b+1.$$
    This would force $a<-1$, which is a contradiction.
    \item If $c\le \max\left\{b, \frac{(a+b)^2}{\,c-1}+1\right\}=\frac{(a+b)^2}{\,c-1}+1$, then $(c-1)^2 \le (a+b)^2.$ This inequality simplifies to $a \ge c-b-1$. 
\end{itemize}
$(ii)$ Suppose, for contradiction, that \(b' \neq b\). Without loss of generality, assume \(b' > b\). 

From \cref{eq:BoundaryOneHoledTorus} we have:
\[
|\tr_{\SL}(\partial T)|
= \frac{2(a + b + c - 1)^2}{(a-1)(b-1)(c-1)} + 2
= \frac{2(a + b' + c - 1)^2}{(a-1)(b' - 1)(c-1)} + 2,
\]
which implies
\[
(b - b') \bigl( (a + c)^2 - (b - 1)(b' - 1) \bigr) = 0.
\]
Since \(b' > b\), it must be that
\[
(a + c)^2 = (b - 1)(b' - 1).
\]
However, note that since \(1 < a \le b < b' \le c\), it follows that
\[
(a+c)^2 > (b'-1)^2 > (b-1)(b'-1),
\]
which leads to a contradiction. The argument is symmetric if \(b' < b\). Thus we conclude \(b' = b\).

\end{proof}

\begin{prop}[Bounds for the orthosystole of a one-holed torus]\label{prop: LowerBoundOsysTorus}
Let $T$ be a one-holed torus of boundary length $L$. Then
$$ \frac{4}{\cosh(L/2)-1}+1  <  \cosh(\operatorname{osys}{T}) \le \frac{1}{\cosh(L/6)- 1} + 1.$$
In addition, if $T$ is an OI one-holed torus, then we have a sharp lower bound
\[
\frac{5}{\cosh(L/2)-1}+1 \leq  \cosh(\operatorname{osys}{T}).
\]
\end{prop}

\begin{proof}
    The upper bound is given by Bavard in \cref{thm:Bavardupperbound}. We will give a proof for the lower bound. Let $(a,b,c),$ where $ 1 < a\le b \le c$, be the cosh-length triple of a minimal orthobasis of $T$. By \cref{lem:propertyofMinimal}, $a=\cosh (\operatorname{osys}(T)) \ge c-b-1$. By \cref{lem:IdentitiesOnOneHoledTorus}, we have that 
    $$2\cosh(L/2)-2 = \frac{2(a+b+c-1)^2}{(a-1)(b-1)(c-1)}.$$ 
It implies that 
$$
    2(\cosh(L/2)-1)(a-1) = \frac{2(a+b+c-1)^2}{(b-1)(c-1)} \ge  \frac{2((c-b-1)+b+c-1)^2}{(b-1)(c-1)}=\frac{8(c-1)}{b-1} \geq 8.
$$
The equality happens if and only if $a =c-b-1$ and $b=c$, which implies that $a=-1$, a contradiction. Rearranging the inequality above, we obtain a strict lower bound 
\begin{equation}\label{equ:InequalityBound}
    \frac{4}{\cosh(L/2)-1}+1 <  a = \cosh(\operatorname{osys}{T}).
\end{equation}

Now assume that $T$ is an OI one-holed torus. By \cref{prop:IntegralityOfdAOISurfaces} and \cref{lem:IdentitiesOnOneHoledTorus}, $\tr_{\SL}^2(\lambda_\alpha)$ is an integer at least $5$ and 
    \begin{equation}\label{eq:lowerboundby10}
\tr_{\SL}^2(\lambda_\alpha)=\frac{(a+b+c-1)^2}{(b-1)(c-1)} = (\cosh(L/2)-1)(a-1)\geq 5,
    \end{equation}
where $\lambda_\alpha$ is the accompanying simple closed geodesic of $\alpha$. This gives us the desired lower bound. Moreover, by applying \cref{eq:BoundaryOneHoledTorus}, one can verify that this bound is attained for OI one-holed tori whose minimal orthobasis cosh-length triples are \((2,13,16)\), \((3,17,21)\), \((6,36,64)\), and \((11,49,61)\).
\end{proof}

Using the ingredients just discussed, we now proceed with the proof of \cref{thm:classificationOItorus}, which provides a complete classification of OI one-holed tori.

\begin{proof}[Proof of \cref{thm:classificationOItorus}]
$(ii) \Rightarrow (i):$ The proof is given in \cite[Theorem 1.1]{doan2021}.\\
$(i) \Rightarrow (iii):$  If $T$ is an OI one-holed torus, then $2 \le a\le b\le c$. Denote by $L$ the boundary length of $T$. By \cref{lem:propertyofMinimal}, $a=\cosh (\operatorname{osys}(T)) \ge c-b-1$. By \cref{prop: LowerBoundOsysTorus}, we have that 
\begin{equation}\label{eq:C_L}
        a\le C_L:=\frac{1}{\cosh(L/6)-1}+1.
    \end{equation}
Thus $c \ge b \ge c-a-1 \ge c-C_L-1$. Combining with \cref{eq:lowerboundby10}, we get 
$$5(c-C_L-2)(c-1) \le 5(b-1)(c-1) \le (a+b+c-1)^2 \le (C_L+b+c-1)^2 \le (C_L+2c-1)^2.$$
Therefore, \begin{equation}\label{equ:c}
    c \le \frac{5C_L+3+\sqrt{29C_L^2+2C_L+53}}{2}.
\end{equation}
By \cref{lem:IdentitiesOnOneHoledTorus}, we have that $|\tr_{\SL}(\partial T)| = 2\cosh(L/2) \in \mathbb{N}$, thus $2\cosh(L/2)\ge 3$, or equivalently,
$$L \ge 2 \operatorname{arcosh}(3/2).$$ 
Combining with \cref{eq:C_L},
$$a \le C_L:=\frac{1}{\cosh(L/6)-1}+1 \le\frac{1}{\cosh\left(\frac{2 \operatorname{arcosh}(3/2)}{6}\right)-1}+1
 \approx 20.2672< 21.$$
Combining with \cref{equ:c}, $$c \le \frac{5C_L+3+\sqrt{29C_L^2+2C_L+53}}{2} < 107.$$
We also have that
$$2\le a \le \frac{1}{\cosh(L/6)-1}+1.$$
which implies $\cosh(L/6) \le 2$ and that 
$$|\tr_{\SL}(\partial T)| = 2\cosh(L/2) = 2(4\cosh^3(L/6)-3\cosh(L/6)) \le 52.$$
Therefore $a,b,c,|\tr_{\SL}(\partial T)| \in \mathbb{N}$ satisfying that $2\le a\le b\le c \le 107$, $c-b-1 \le a \le 20$ and $3 \le |\tr_{\SL}(\partial T)| \le 52.$ \\
$(iii) \Rightarrow (ii):$ One can verify using a short computer program that the list provided in $(ii)$ contains all possible solutions satisfying condition~$(iii)$.
\end{proof}

\begin{remark}
Based on the classification of OI one-holed tori, \cref{table:1} presents a list of the associated Markoff-type equations, where the initial solution corresponds to the traces of the triple of simple closed geodesics that accompany a minimal orthobasis of the one-holed torus.
\end{remark}
 \begin{longtable}{|c|c|c|}
\hline
% \cellcolor{gray!30}
Orthobasis  & Markoff-type equation & Initial solution\\
% \hline
\endfirsthead

% \hline
% \cellcolor{gray!30} \textbf{Orthobasis}  & \cellcolor{gray!30} \textbf{Markoff-type equation} & \cellcolor{gray!30} \textbf{Initial solution}\\
% \hline
\endhead

% \hline
% \endfoot

% \hline
% \endlastfoot

  \hline
  $(2,2,2)$ & $X^2+Y^2+Z^2-XYZ=-50$ & $(5,5,5)$ \\
  \hline
  $(2,2,3)$ & $X^2+Y^2+Z^2-XYZ=-36$ & $(3\sqrt{2},3\sqrt{2},6)$\\
  \hline
  $(2,2,5)$ & $X^2+Y^2+Z^2-XYZ=-32$& $(4,4,8)$\\
  \hline
  $(2,3,6)$ & $X^2+Y^2+Z^2-XYZ=-20$& $(\sqrt{10},2\sqrt{5},5\sqrt{2})$\\
  \hline
  $(2,4,4)$ & $X^2+Y^2+Z^2-XYZ=-18$& $(3,3\sqrt{3},3\sqrt{3})$\\
  \hline
  $(2,4,7)$ & $X^2+Y^2+Z^2-XYZ=-16$& $(2\sqrt{2},2\sqrt{6},4\sqrt{3})$\\
  \hline
  $(2,5,8)$ & $X^2+Y^2+Z^2-XYZ=-14$& $(\sqrt{7},2\sqrt{7},7)$\\
  \hline
  $(2,7,10)$ & $X^2+Y^2+Z^2-XYZ=-12$& $(\sqrt{6},6,3\sqrt{6})$\\
  \hline
  $(2,13,16)$ & $X^2+Y^2+Z^2-XYZ=-10$& $(\sqrt{5},2\sqrt{15},5\sqrt{3})$\\

  \hline
  $(3,3,3)$ & $X^2+Y^2+Z^2-XYZ=-16$& $(4,4,4)$\\
  \hline
  $(3,3,7)$ & $X^2+Y^2+Z^2-XYZ=-12$& $(2\sqrt{3},2\sqrt{3},6)$\\
  \hline
  $(3,5,5)$ & $X^2+Y^2+Z^2-XYZ=-9$& $(3,3\sqrt{2},3\sqrt{2})$\\
  \hline
  $(3,5,9)$ & $X^2+Y^2+Z^2-XYZ=-8$& $(2\sqrt{2},4,4\sqrt{2})$\\
  \hline
  $(3,9,13)$ & $X^2+Y^2+Z^2-XYZ=-6$ & $(\sqrt{6},2\sqrt{6},6)$\\
  \hline
  $(3,17,21)$ & $X^2+Y^2+Z^2-XYZ=-5$& $(\sqrt{5},2\sqrt{10},5\sqrt{2})$\\
  \hline
  $(4,4,5)$ & $X^2+Y^2+Z^2-XYZ=-8$& $(2\sqrt{3},2\sqrt{3},4)$\\
  \hline
  $(4,5,10)$ & $X^2+Y^2+Z^2-XYZ=-6$& $(3,2\sqrt{3},3\sqrt{3})$\\
  \hline
  $(4,6,6)$ & $X^2+Y^2+Z^2-XYZ=-6$& $(3,\sqrt{15},\sqrt{15})$\\
  \hline
  $(4,11,16)$ & $X^2+Y^2+Z^2-XYZ=-4$& $(\sqrt{6},2\sqrt{5},\sqrt{30})$\\
  \hline
  $(5,5,11)$ & $X^2+Y^2+Z^2-XYZ=-5$& $(\sqrt{10},\sqrt{10},5)$\\
  \hline
  $(5,7,13)$ & $X^2+Y^2+Z^2-XYZ=-4$& $(2\sqrt{2},2\sqrt{3},2\sqrt{6})$\\
  \hline
  $(5,13,19)$ & $X^2+Y^2+Z^2-XYZ=-3$& $(\sqrt{6},3\sqrt{2},3\sqrt{3})$\\
  \hline
  $(6,6,9)$ & $X^2+Y^2+Z^2-XYZ=-4$& $(\sqrt{10},\sqrt{10},4)$\\
  \hline
  $(6,36,64)$ & $X^2+Y^2+Z^2-XYZ=-2$& $(\sqrt{5},\sqrt{35},3\sqrt{7})$\\
  \hline
  $(7,9,9)$ & $X^2+Y^2+Z^2-XYZ=-3$& $(3,2\sqrt{3},2\sqrt{3})$\\
  \hline
  $(7,17,25)$ & $X^2+Y^2+Z^2-XYZ=-2$& $(\sqrt{6},4,2\sqrt{6})$\\
  \hline
  $(8,13,22)$ & $X^2+Y^2+Z^2-XYZ=-2$& $(\sqrt{7},2\sqrt{3},\sqrt{21})$\\
  \hline
  $(9,11,21)$ & $X^2+Y^2+Z^2-XYZ=-2$& $(2\sqrt{2},\sqrt{10},2\sqrt{5})$\\
  \hline
  $(10,10,17)$ & $X^2+Y^2+Z^2-XYZ=-2$& $(3,3,4)$\\
  \hline
  $(10,12,12)$ & $X^2+Y^2+Z^2-XYZ=-2$& $(3,\sqrt{11},\sqrt{11})$\\
  \hline
  $(11,49,61)$ & $X^2+Y^2+Z^2-XYZ=-1$& $(\sqrt{5},2\sqrt{6},\sqrt{30})$\\
  \hline
  $(13,29,43)$ & $X^2+Y^2+Z^2-XYZ=-1$& $(\sqrt{6},\sqrt{14},\sqrt{21})$\\
  \hline
  $(17,19,37)$ & $X^2+Y^2+Z^2-XYZ=-1$& $(2\sqrt{2},3,3\sqrt{2})$\\
  \hline
   $(19,21,21)$ & $X^2+Y^2+Z^2-XYZ=-1$& $(3,\sqrt{10},\sqrt{10})$\\
  \hline
\caption{Markoff-type equations of OI one-holed tori.}
\label{table:1} 
\end{longtable}

\section{Ortho-integral surfaces and quadratic forms}\label{sec:Ortho-integral surfaces and quadratic forms}

In this section, we describe some connections between OI one-holed tori and OI pairs of pants with quadratic forms. We leverage these connections to construct examples of infinitely many pairwise non-commensurable OI surfaces and to compute arithmetic invariants of a one-holed torus and a pair of pants, from the data of a decomposition of the surface into isometric right-angled hexagons.  

\subsection{Connection with quadratic form}
Let $H(a,b,c)$ be a right-angled hexagon in $\mathbb{H}^2$ with three alternating sides of cosh-length $a$, $b$ and $c$ in counterclockwise order. Since cosh of any positive number exceeds 1, we have $a,b,c >1$. Associated with this right-angled hexagon, we consider the quadratic form $q = q_{(a,b,c)}$ given by
\begin{equation}
\label{eq:QuadraticForm}
q_{(a,b,c)}(x,y,z) = (a^2-1)x^2 + (b^2-1)y^2 + (c^2-1)z^2 - 2(ab+c)xy-2(bc+a)yz -2(ca+b)zx, 
\end{equation}
where $a,b,c > 1$. The bilinear form associated to $q$ is 
\begin{equation}
\label{eq:BilinearForm}
B({\bf x},{\bf y}) ={\bf x}^\top Q {\bf y} = {\bf x}^\top \begin{pmatrix}
    a^2-1 & -ab-c & -ca-b\\
    -ab-c & b^2-1 & -bc-a\\
    -ca-b & -bc-a & c^2-1
\end{pmatrix}{\bf y}.
\end{equation}
Also associated with this right-angled hexagon, we consider the following collection of vectors in $\mathbb{R}^3$.
\begin{equation}
\label{eq:SpecialVectors}
\begin{aligned}
    {\bf e}_1 &= (1,0,0)^\top,   &&{\bf e}_2 = (0,1,0)^\top,  &&&{\bf e}_3 = (0,0,1)^\top \\
    {\bf u}_\alpha &= (-1,c,b)^\top,  && {\bf u}_\beta = (c,-1,a)^\top, &&&{\bf u}_\gamma =  (b,a,-1)^\top \\
    {\bf v}_\alpha &= (0,c^2-1,bc+a)^\top,  && {\bf v}_\beta = (ca+b,0,a^2-1)^\top, &&&{\bf v}_\gamma =  (b^2-1,ab+c,0)^\top \\
    {\bf w}_\alpha &= (0,bc+a,b^2-1)^\top,  && {\bf w}_\beta = (c^2-1,0,ca+b)^\top, &&&{\bf w}_\gamma =  (ab+c,a^2-1,0)^\top \\
\end{aligned}   
\end{equation} 
\noindent
A straight forward calculation shows that each of the following set is an orthogonal basis for $(V,q)$.
\begin{equation}
\label{eq:OrthogonalRelations}
    \{{\bf v}_\beta,{\bf u}_\beta,{\bf e}_1\},  \{{\bf w}_\gamma,{\bf u}_\gamma,{\bf e}_1\}, 
    \{{\bf v}_\gamma,{\bf u}_\gamma,{\bf e}_2\},  \{{\bf w}_\alpha,{\bf u}_\alpha,{\bf e}_2\}, 
    \{{\bf v}_\alpha,{\bf u}_\alpha,{\bf e}_3\},  \{{\bf w}_\beta,{\bf u}_\beta,{\bf e}_3\}  
\end{equation}
where $q({\bf u}_*) > 0 $, $q({\bf v}_*) < 0 $, $q({\bf w}_*) < 0$  and $q({\bf e}_i) > 0 $ for all $* \in \{\alpha,\beta,\gamma\}$ and $1 \leq i \leq 3$. Therefore, the quadratic form $q_{(a,b,c)}$ has signature $(+,+,-)$ for all real numbers $a,b,c > 1 $. Since $q({\bf v}_\alpha) <0$, we can choose a model for the hyperbolic plane by projectivizing the component of $q({\bf x})<0$ that contains ${\bf v}_\alpha$. The hyperbolic plane is modeled on the projectivization $P(\mathcal{C}^+)$ of the set 
\[
\mathcal{C}^+ = \{ {\bf x } = (x_1,x_2,x_3)^\top \mid q({\bf x}) < 0 \text{ and } (c^2-1)x_2+ (bc+a)x_3 >0  \}
\]
where the hyperbolic distance on $P(\mathcal{C}^+)$ is given by
\begin{equation}
    \label{eq:DistanceH2}
    d_{\mathbb{H}^2}({\bf x },{\bf y }) = \operatorname{arcosh}\left(\frac{-B({\bf x },{\bf y })}{\sqrt{B({\bf x },{\bf x })B({\bf y },{\bf y })}}\right).
\end{equation}
To simplify notation, we identify a point in $P(\mathcal{C}^+)$ to a vector in $\mathbb{R}^3$ that projects onto it. The following proposition gives some geometric interpretation of vectors considered in \cref{eq:SpecialVectors}. 

\begin{prop}
\label{prop:GeometricInterpretationOfuvw}
The vectors ${\bf v}_*$ and ${\bf w}_*$ project to six points on $P(\mathcal{C}^+)$ whose convex hull is an all-right angled hexagon isometric to $H(a,b,c)$. The geodesic segments $\alpha = [{\bf v}_\beta,{\bf w}_\gamma]$, $\beta= [{\bf v}_\gamma,{\bf w}_\alpha]$, and $\gamma=[{\bf v}_\alpha,{\bf w}_\beta]$ has cosh-length
\[
\cosh(\ell(\alpha)) = a, \ \cosh(\ell(\beta)) = b,  
\cosh(\ell(\gamma)) = c.
\]
The orthogonal complement ${\bf u}_*^\perp \subset \mathbb{R}^3$ projects to an (infinite) geodesic in  $P(\mathcal{C}^+)$ that contains the geodesic segment $[{\bf v}_*, {\bf w}_*]$ for each $* \in \{\alpha,\beta,\gamma\}$. The orthogonal complement ${\bf e}_i^\perp \subset \mathbb{R}^3$ projects to an (infinite) geodesic in  $P(\mathcal{C}^+)$ that contains the geodesic segment $*$ for each $* \in \{\alpha,\beta,\gamma\}$.
\end{prop}

\begin{proof}
Since $q({\bf v}_*) < 0 $ and $q({\bf w}_*) < 0$ for any $* \in \{\alpha,\beta,\gamma\}$, the vectors ${\bf v}_*$ and ${\bf w}_*$ project to six points on $P(\mathcal{C}^+)$. Since $q({\bf u}_*) > 0 $ and $q({\bf e}_i) > 0$ for any $* \in \{\alpha,\beta,\gamma\}$ and $1\leq i \leq 3$, the orthogonal complement ${\bf u}_*^\perp$ and ${\bf e}_i^\perp$ projects to geodesics on $P(\mathcal{C}^+)$. The fact that each set in \cref{eq:OrthogonalRelations} is an orthogonal basis for $(V,q)$ implies that the convex hull of $\{{\bf v}_\alpha,{\bf v}_\beta,{\bf v}_\gamma,{\bf w}_\alpha,{\bf w}_\beta,{\bf w}_\gamma\}$ is an all-right angled hexagon and that the orthogonal complements ${\bf u}_*^\perp$ and ${\bf e}_i^\perp$ projects to the sides of the hexagon as claimed.  Using \cref{eq:DistanceH2}, we compute the cosh of the distance between the points ${\bf v}_*$ and ${\bf w}_*$ to be
\[
\cosh(d_{\mathbb{H}^2}({\bf v}_\beta,{\bf w}_\gamma)) =  a, \ \cosh(d_{\mathbb{H}^2}({\bf v}_\gamma,{\bf w}_\alpha)) = b, \ \cosh(d_{\mathbb{H}^2}({\bf v}_\alpha,{\bf w}_\beta)) = c.
\] which implies that the convex hull of $\{{\bf v}_\alpha,{\bf v}_\beta,{\bf v}_\gamma,{\bf w}_\alpha,{\bf w}_\beta,{\bf w}_\gamma\}$ in $P(\mathcal{C}^+)$ is isometric to an all-right angled hexagon $H(a,b,c)$. 
\end{proof}

\begin{figure}[h!]
    \centering
\begin{tikzpicture}[scale = 2]
    \draw[thick] (0,0) circle (1);
    \draw[thick] (-0.25,{sqrt(3)/4}) arc (240:300:0.5);
    \draw[thick] (-0.25,{sqrt(3)/4}) arc (0:-60:0.5);
    \draw[thick] (-0.25,{-sqrt(3)/4}) arc (0:60:0.5);
    \draw[thick] (-0.25,{-sqrt(3)/4}) arc (120:60:0.5);
    \draw[thick] (0.5,0) arc (240:180:0.5);
    \draw[thick] (0.5,0) arc (120:180:0.5);
    
    \filldraw [black] (0,{(sqrt(3)+1)/2}) circle (0.5pt) node[anchor=south]{${\bf e}_3$};
    \filldraw [black] ({-(sqrt(3)+3)/4},{-(sqrt(3)+1)/4}) circle (0.5pt) node[anchor=east]{${\bf e}_1$};
    \filldraw [black] ({(sqrt(3)+3)/4},{-(sqrt(3)+1)/4}) circle (0.5pt) node[anchor=west]{${\bf e}_2$};
    \filldraw [black] ({-(sqrt(3)+3)/4},{(sqrt(3)+1)/4}) circle (0.5pt) node[anchor=east]{${\bf u}_\beta$};
    \filldraw [black] ({(sqrt(3)+3)/4},{(sqrt(3)+1)/4}) circle (0.5pt) node[anchor=west]{${\bf u}_\alpha$};
    \filldraw [black] (0,-{(sqrt(3)+1)/2}) circle (0.5pt) node[anchor=north]{${\bf u}_\gamma$};
    %%%%%%%%%%%%%%%%%%%%%
    \filldraw [black] (0.5,0) circle (0.5pt) node[anchor=west]{${\bf w}_\alpha$};
    \filldraw [black] (0.25,{sqrt(3)/4}) circle (0.5pt) node[anchor=west]{${\bf v}_\alpha$};
    \filldraw [black] (-0.25,{sqrt(3)/4}) circle (0.5pt) node[anchor=east]{${\bf w}_\beta$};
    \filldraw [black] (-0.5,0) circle (0.5pt) node[anchor=east]{${\bf v}_\beta$};
    \filldraw [black] (-0.25,{-sqrt(3)/4}) circle (0.5pt) node[anchor=east]{${\bf w}_\gamma$};
    \filldraw [black] (0.25,{-sqrt(3)/4}) circle (0.5pt) node[anchor=west]{${\bf v}_\gamma$};
    %%%%%%%%%%%%%%%%%%%%
\end{tikzpicture}  
    \caption{The all right-angled hexagon $H(a,b,c)$. }
    \label{fig:RightAngledHexagon}
\end{figure}
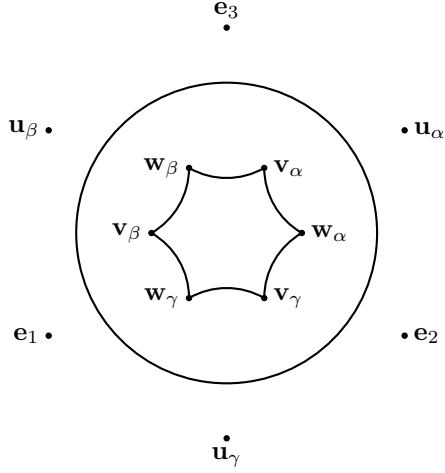

We identify $H(a,b,c)$ with the convex hull of $\{{\bf v}_\alpha,{\bf v}_\beta,{\bf v}_\gamma,{\bf w}_\alpha,{\bf w}_\beta,{\bf w}_\gamma\}$ in $P(\mathcal{C}^+)$. We summarize \cref{prop:GeometricInterpretationOfuvw} in \cref{fig:RightAngledHexagon}. 

Naturally, for $* \in \{\alpha,\beta,\gamma\}$ the matrices $f_*$ and  $g_*$ from \cref{thm:GroupAssociatedToOrthoGeodesics} also have a geometric interpretation. In particular, $f_*$ and $g_*$ preserve the quadratic form $q$ and $\mathcal{C}^+$, thus defining non-trivial isometries of the hyperboloid model. Since $f_\alpha$ fixes ${\bf v}_\beta$, ${\bf w}_\gamma$, the isometry $f_\alpha$ is a reflection across $\alpha$. Similarly, $f_\beta$ is a reflection across $\beta$, and $f_\gamma$ is a reflection across $\gamma$. The matrix $g_*$ has order 2 and $\det(g_*) = 1$ for any $* \in \{\alpha,\beta,\gamma \}$, and thus defines an order two rotation. Since $g_\alpha({\bf v}_\beta) = {\bf w}_\gamma$, $g_\alpha$ is an order two rotation about the midpoint of $\alpha$. Similarly, $g_\beta$ and $g_\gamma$ are order two rotations about the midpoints of $\beta$ and $\gamma$ respectively. The following proposition is an immediate consequence of the discussion above. The key idea of \cref{prop:PantsAndOneHoledTorusFromRootFlipping} was communicated to us by Greg McShane. We would like to thank him for sharing his insights.

\begin{prop}
\label{prop:PantsAndOneHoledTorusFromRootFlipping}
    Let $P(a,b,c)$ be a pair of pants with a standard orthobasis $\mathcal{H}_P = \{\alpha,\beta,\gamma\}$ and $(a,b,c):= (\cosh\ell(\alpha),\cosh\ell(\beta),\cosh\ell(\gamma))$. The group $G_P$ generated by $f_\alpha$, $f_\beta$ and $f_\gamma$ is a discrete subgroup of $\ORTH^+(q_{a,b,c};\mathbb{R})$. The universal cover $\widetilde{P}(a,b,c)$ of $P(a,b,c)$ is isometric to a tiling of a convex subset of $P(\mathcal{C}^+)$ by the orbit of $H(a,b,c)$ under the action of $G_P$. The subgroup $\Gamma_P := \langle f_\alpha f_\beta, f_\alpha f_\gamma \rangle$ has index two in $G_P$, and $\widetilde{P}(a,b,c)/\Gamma_P$ is isometric to $P(a,b,c)$.
    
    Let $T(a,b,c)$ be a one-holed torus with a standard orthobasis $\mathcal{H}_T = \{\alpha,\beta,\gamma\}$ and $(a,b,c):= (\cosh\ell(\alpha),\cosh\ell(\beta),\cosh\ell(\gamma))$. The group $G_T$ generated by $g_\alpha$, $g_\beta$ and $g_\gamma$ is a discrete subgroup of $\ORTH^+(q_{a,b,c};\mathbb{R})$. The universal cover $\widetilde{T}(a,b,c)$ of $T(a,b,c)$ is isometric to a tiling of a convex subset of $P(\mathcal{C}^+)$ by the orbit of $H(a,b,c)$ under the action of $G_T$. The subgroup $\Gamma_T := \langle g_\alpha g_\beta, g_\alpha g_\gamma \rangle$ has index two in $G_T$, and $\widetilde{T}(a,b,c)/\Gamma_T$ is isometric to $T(a,b,c)$.
\end{prop}

\subsection{Verifying arithmeticity of the double of OI pair of pants and OI one-holed torus}

In this section, we compute the arithmetic invariants of OI one-holed torus and OI pair of pants. We verify that the doubles of OI one-holed tori and OI pairs of pants are all arithmetic surfaces of genus 2. The following theorem provides a criterion for determining whether a finite-covolume Fuchsian group is arithmetic.

\begin{thm}[\cite{takeuchi1975characterization}]
\label{thm:TakeuchiArithmeticCriterion}
    Let $\Gamma$ be a finite-covolume Fuchsian group. Then $\Gamma$ is arithmetic if and only if the following three conditions hold
    \begin{itemize}
        \item The invariant trace field $k\Gamma$ is a totally real field.
        \item $\tr(\gamma)$ is an algebraic integer for every $\gamma \in \Gamma$.
        \item $A\Gamma$ is ramified at all real places of $k\Gamma$ except one. 
    \end{itemize}
\end{thm}

Here, the \emph{invariant trace field $k\Gamma$} is $k\Gamma = \mathbb{Q}(\{ \tr_{\SL} \gamma^2 \mid \gamma \in \Gamma \})$. The \emph{invariant quaternion algebra $A\Gamma$} is defined to be $A\Gamma = \{\sum a_i \gamma_i \mid a_i \in k\Gamma, \gamma_i \in \Gamma\}$. To compute the arithmetic invariants of the genus-two surface obtained by doubling either a pair of pants $P$ or a one-holed torus $T$, we need the reflections across the boundary components of $P$ and $T$. To that end, we consider the following matrices
\begin{equation}
    \label{eq:ReflectionsAcrossBoundary}
    r_\alpha = \begin{pmatrix} -1 & 0 & 0 \\ 2c & 1 & 0 \\ 2b & 0 & 1 \end{pmatrix}, \ 
    r_\beta = \begin{pmatrix} 1 & 2c & 0 \\ 0 & -1 & 0 \\ 0 & 2a & 1 \end{pmatrix}, \
    r_\gamma = \begin{pmatrix} 1 & 0 & 2b \\ 0 & 1 & 2a \\ 0 & 0 & -1 \end{pmatrix}.
\end{equation}
For $*\in\{\alpha,\beta,\gamma\}$, we can check that $r_*$ fixes ${\bf v}_*$ and ${\bf w}_*$ and defines reflections across the sides of $H(a,b,c)$ opposite to $*$, see \cref{fig:RightAngledHexagon}.

The description of OI one-holed tori and OI pairs of pants in terms of an orthobasis makes it natural to describe the corresponding Fuchsian groups using the quadratic form $q=q_{(a,b,c)}$. 

\begin{prop}
\label{prop:DoubleOfPantsAndTorus}
    Consider the following groups 
    \begin{equation}
    \label{eq:DoubleOfPantsAndTorus}
        \Delta_P = \langle f_\alpha f_\beta, f_\alpha f_\gamma, r_\alpha r_\beta, r_\alpha  r_\gamma\rangle, \text{ and } \Delta_T = \langle g_\alpha g_\beta,  g_\alpha g_\gamma, r_\alpha g_\alpha g_\beta r_\alpha, r_\alpha g_\alpha g_\gamma r_\alpha  \rangle.
    \end{equation}    Then $\mathbb{H}^2 /\Delta_P$ and $\mathbb{H}^2 /\Delta_T$ are hyperbolic surfaces of genus 2 isometric to the double of the pair of pants $P(a,b,c)$ and the one-holed torus $T(a,b,c)$, respectively.
\end{prop}

We recall the relationship between the trace of matrices in $\PSL_2(\mathbb{R})$ and the trace of matrices in $\SO^+(2,1;\mathbb{R})$. The Lie algebra $\liesl_2(\mathbb{R})$ is a 3-dimensional real vector space of traceless $2\times 2$ matrices. The algebra $\liesl_2(\mathbb{R})$ is equipped with a bilinear form of signature $(2,1)$ given by $(X,Y) \mapsto \tr(XY)$. The Lie group $\PSL_2(\mathbb{R})$ acts on the Lie algebra $\liesl_2(\mathbb{R})$ by conjugation respecting the bilinear form. Therefore, the action of conjugation of the group $\PSL_2(\mathbb{R})$ defines an isomorphism with $\SO^+(2,1;\mathbb{R})$. With respect to the basis 
\[
\begin{pmatrix}
    1 & 0 \\ 0 & -1
\end{pmatrix}, \ 
\begin{pmatrix}
    0 & 1 \\ 1 & 0
\end{pmatrix}, \ 
\begin{pmatrix}
    0 & 1 \\ -1 & 0
\end{pmatrix}
\]
of $\liesl_2(\mathbb{R})$, the isomorphism between $\PSL_2(\mathbb{R})$ and $\SO^+(2,1;\mathbb{R})$ is given by 
\[
\Psi\left[\pm \begin{pmatrix} a & b \\ c & d \end{pmatrix}\right] = 
\begin{pmatrix}
    bc + ad & -ac + bd & -ac - bd \\
    cd - ab & (a^2-b^2-c^2+d^2)/2 & (a^2 + b^2 - c^2 - d^2)/2 \\
   -cd - ab & (a^2-b^2+c^2-d^2)/2 & (a^2+b^2+c^2 + d^2)/2
\end{pmatrix}.
\]
For any matrix $A \in \PSL_2(\mathbb{R})$, we have
\begin{equation}
    \label{eq:TraceBetweenSOandSL}
    \tr_{\SO}(\Psi(A)) = a^2 + d^2 + ad + bc = \tr_{\SL}(A)^2 - 1.
\end{equation}
Since trace is a conjugation invariant, \cref{eq:TraceBetweenSOandSL} holds independent of the choice of quadratic form for $\mathbb{R}^3$. 

To compute the arithmetic invariant, suppose that $a,b,c$ are integers greater than 1. Then $(V,q_{(a,b,c)})$ is a quadratic space defined over $\mathbb{Q}$ of signature $(2,1)$. We identify $V$ with $\mathbb{Q}^3$ by choosing a basis such that the matrix associated to $q$ is 
\[
\begin{pmatrix}
    a^2-1 & -ab-c & -ac-b\\
    -ab-c & b^2-1 & -bc-a\\
    -ac-b & -bc-a & c^2-1
\end{pmatrix}
\]
The following vectors in $V$ give an orthonormal basis for $(V,q_{(a,b,c)})$
\begin{equation}
    \mathbf{w}_\alpha= \begin{pmatrix} 0 \\ bc+a \\ b^2-1 \end{pmatrix}, \quad 
    \mathbf{u}_\alpha = \begin{pmatrix}
        -1 \\ c \\ b
    \end{pmatrix},  \quad \text{and} \quad
    \mathbf{e}_2 = \begin{pmatrix}
        0 \\ 1 \\ 0
    \end{pmatrix}. 
\end{equation}

By \cref{eq:QuaternionAlgebraAssociatedToQuadraticSpace}, the quaternion algebra associated to the quadratic space $(V,q_{(a,b,c)})$ is 
\begin{equation}\label{equ:ComputingIQA}
\left(\frac{-q(\mathbf{w}_\alpha)q(\mathbf{u}_\alpha),-q(\mathbf{w}_\alpha)q(\mathbf{e}_2)}{\mathbb{Q}}\right) \cong \left(\frac{(b^2-1),(a^2+b^2+c^2+2abc-1)}{\mathbb{Q}}\right)
\end{equation}

\begin{prop} \label{prop:Double}
The double of any OI pair of pants or any OI one-holed torus is an arithmetic surface of genus two.
\end{prop}

\begin{proof}
    The proof is a direct computation verifying the three conditions in \cref{thm:TakeuchiArithmeticCriterion}. Let $P$ be an OI pair of pants and $T$ be an OI one-holed torus. The Fuchsian group that corresponds to the double of $P$ and the double of $T$ is $\Delta_P$ and $\Delta_T$ given by \cref{prop:DoubleOfPantsAndTorus}. Using \cref{eq:TraceBetweenSOandSL}, \cref{eq:DoubleOfPantsAndTorus} and the classification of OI pairs of pants and OI one-holed tori, we verify that the trace set of $\Delta_P$ and $\Delta_T$ is integral. After doubling $P$ or $T$, every orthogeodesics in $P$ and in $T$ becomes a closed geodesic in the genus two surface, the trace of these particular geodesics is twice the cosh of the length of the orthogeodesic from which they originate. To check that the trace set of $\Delta_P$ and of $\Delta_T$ is integral, we use \cite[Lemma 3.5.2]{maclachlan2003book} which produces a finite set of traces that we need to check for integrality. In particular, for a 4-generated group $\Gamma = \langle w,x,y,z 
    \rangle$, it suffices to check integrality for the traces of the elements
    \[
    \{w,x,y,z,wx,wy,wz,xy,xz,yz,wxy,wxz,wyz,xyz,wxyz\}.
    \]
    In all cases, the traces of all elements in $\Delta_P$ or $\Delta_T$ are all algebraic integers. See \cref{tab:OIPantsInvariants} and \cref{tab:OIToriInvariants} for the resulting computations. 
    
    To compute the invariant trace field, \cite[Lemma 3.5.3]{maclachlan2003book} and \cite[Lemma 3.5.5]{maclachlan2003book} imply that for a non-elementary subgroup $\Gamma = \langle w,x,y,z 
    \rangle \leq \PSL_2(\mathbb{C})$, it suffices to compute the traces of   
    \[
    \{w^2,x^2,y^2,z^2,w^2x^2,w^2y^2,w^2z^2,x^2y^2,x^2z^2,y^2z^2,w^2x^2y^2,w^2x^2z^2,w^2y^2z^2,x^2y^2z^2\}.
    \]
    The resulting computations show that the invariant trace field of $\Delta_P$ and $\Delta_T$ are all $\mathbb{Q}$ which are totally real. The last condition in \cref{thm:TakeuchiArithmeticCriterion} is vacuous. 
\end{proof}

The following tables record the arithmetic data of OI pair of pants and one-holed torus. The first column records the cosh-length of an ortho-basis for the surfaces. The second column records a Hilbert symbol of the invariant quaternion algebra of the surface computed using \cref{equ:ComputingIQA}. The third column records the finite primes in $\mathbb{Z}$ over which $A$ ramifies, i.e. the completion of $A$ at $p$ is a division algebra. For a fixed number field, the ramification set is a complete invariant of the invariant quaternion algebra \cite[Theorem 2.7.5]{maclachlan2003book} and therefore uniquely specifies the algebra up to isomorphism. The final column records the trace set of the genus 2 surface group that corresponds to doubling $P$ or $T$.

\begin{longtable}{|c|c|c|c|}
\hline
Orthobasis & \(A/\mathbb{Q}\) & \(\operatorname{Ram}(A)\) & \(\mathrm{tr}(\Delta_P)\) \\
\hline
\endfirsthead

% \hline
% \textbf{Orthobasis} & \(\boldsymbol{A/\mathbb{Q}}\) & \(\boldsymbol{Ram(A)}\) & \(\boldsymbol{\mathrm{tr}(D\Gamma)}\) \\
% \hline
\endhead

% \hline
% \endfoot

\hline
% \endlastfoot
         $(2,2,2)$ & $(27,3)=(3,3)$ & $\{2,3\}$ & $\mathbb{Z}$ \\  \hline
         $(2,2,5)$ & $(72,24) = (2,6)$ & $\{2,3\}$ & $\mathbb{Z}[\sqrt{2}]$  \\ \hline
         $(2,2,17)$ & $(432,3288) =(3,2)$& $\{2,3\}$ & $\mathbb{Z}[\sqrt{6}]$ \\ \hline

         $(3,3,3)$ & $(80,8) = (5,2)$ & $\{2,5\}$ & $\mathbb{Z}$\\ \hline
         $(3,3,7)$ & $(192,48) = (3,3) $ & $\{2,3\}$ & $\mathbb{Z}[\sqrt{6}]$\\ \hline
         $(3,3,19)$ & $(720,360) = (10,5)$ & $\{2,5\}$ & $\mathbb{Z}[\sqrt{5}]$\\ \hline
         $(5,5,11)$ & $(720,120) = (5,30)$& $\{2,3\}$ & $\mathbb{Z}[\sqrt{5}]$ \\
\hline 
     \caption{Classification of OI pairs of pants with associated invariants}
    \label{tab:OIPantsInvariants}
\end{longtable}

\begin{longtable}{|c|c|c|c|}

\hline
Orthobasis & \(A/\mathbb{Q}\) & \(\operatorname{Ram}(A)\) & \(\mathrm{tr}(\Delta_T)\) \\
\hline
\endfirsthead

% \hline
% \textbf{Orthobasis} & \(\boldsymbol{A/\mathbb{Q}}\) & \(\boldsymbol{Ram(A)}\) & \(\boldsymbol{\mathrm{tr}(D\Gamma)}\) \\
% \hline
\endhead

% \hline
% \endfoot

\hline
% \endlastfoot

\((2,2,2)\) & \((27,3)=(3,3)\) & \(\{2,3\}\) & \(\mathbb{Z}\) \\
\hline
\((2,2,3)\) & \((40,8) = (10,2)\) & \(\{2,5\}\) & \(\mathbb{Z}[\sqrt{2}]\) \\
\hline
\((2,2,5)\) & \((72,24) = (2,6)\) & \(\{2,3\}\) & \(\mathbb{Z}\) \\
\hline
\((2,3,6)\) & \((120,35) =(30,35)\) & \(\{3,5\}\) & \(\mathbb{Z}[\sqrt{2},\sqrt{5}]\) \\
\hline
\((2,4,4)\) & \((99,15) =(11,15)\) & \(\{2,3\}\) & \(\mathbb{Z}[\sqrt{3}]\) \\
\hline
\((2,4,7)\) & \((180,48) = (5,3)\) & \(\{3,5\}\) & \(\mathbb{Z}[\sqrt{2},\sqrt{3}]\) \\
\hline
\((2,5,8)\) & \((252,63) = (7,7)\) & \(\{2,7\}\) & \(\mathbb{Z}[\sqrt{7}]\) \\
\hline
\((2,7,10)\) & \((432,99) = (3,11)\) & \(\{2,3\}\) & \(\mathbb{Z}[\sqrt{6}]\) \\
\hline
\((2,13,16)\) & \((1260,255)=(35,255)\) & \(\{2,3,5,7\}\) & \(\mathbb{Z}[\sqrt{3},\sqrt{5}]\) \\
\hline
\((3,3,3)\) & \((80,8) = (5,2)\) & \(\{2,5\}\) & \(\mathbb{Z}\) \\
\hline
\((3,3,7)\) & \((192,48) = (3,3)\) & \(\{2,3\}\) & \(\mathbb{Z}[\sqrt{3}]\) \\
\hline
\((3,5,5)\) & \((208,24) = (13,6)\) & \(\{2,13\}\) & \(\mathbb{Z}[\sqrt{2}]\) \\
\hline
\((3,5,9)\) & \((384,80)=(6,5)\) & \(\{2,3\}\) & \(\mathbb{Z}[\sqrt{2}]\) \\
\hline
\((3,9,13)\) & \((960,168) = (15,42)\) & \(\{3,5\}\) & \(\mathbb{Z}[\sqrt{6}]\) \\
\hline
\((3,17,21)\) & \((2880,440) = (5,110)\) & \(\{2,5\}\) & \(\mathbb{Z}[\sqrt{2},\sqrt{5}]\) \\
\hline
\((4,4,5)\) & \((216,24) = (6,6)\) & \(\{2,3\}\) & \(\mathbb{Z}[\sqrt{3}]\) \\
\hline
\((4,5,10)\) & \((540,99) = (15,11)\) & \(\{2,3\}\) & \(\mathbb{Z}[\sqrt{3}]\) \\
\hline
\((4,6,6)\) & \((375,35) = (15,35)\) & \(\{2,3\}\) & \(\mathbb{Z}[\sqrt{15}]\) \\
\hline
\((4,11,16)\) & \((1800,255)= (2,255)\) & \(\{3,5\}\) & \(\mathbb{Z}[\sqrt{5},\sqrt{6}]\) \\
\hline
\((5,5,11)\) & \((720,120) = (5,30)\) & \(\{2,3\}\) & \(\mathbb{Z}[\sqrt{10}]\) \\
\hline
\((5,7,13)\) & \((1152,168)=(2,42)\) & \(\{2,3\}\) & \(\mathbb{Z}[\sqrt{2}, \sqrt{3}]\) \\
\hline
\((5,13,19)\) & \((3024,360)=(21,10)\) & \(\{2,7\}\) & \(\mathbb{Z}[\sqrt{2}, \sqrt{3}]\) \\
\hline
\((6,6,9)\) & \((800,80)=(2,5)\) & \(\{2,5\}\) & \(\mathbb{Z}[\sqrt{10}]\) \\
\hline
\((6,29,36)\) & \((14700,1295) = (3,1295)\) & \(\{2,3,5,7\}\) & \(\mathbb{Z}[\sqrt{5},\sqrt{7}]\) \\
\hline
\((7,9,9)\) & \((1344,80) = (84,5)\) & \(\{3,7\}\) & \(\mathbb{Z}[\sqrt{3}]\) \\
\hline
\((7,17,25)\) & \((6912,624)=(3,39)\) & \(\{2,3\}\) & \(\mathbb{Z}[\sqrt{6}]\) \\
\hline
\((8,13,22)\) & \((5292,483)=(3,483)\) & \(\{2,7\}\) & \(\mathbb{Z}[\sqrt{3},\sqrt{7}]\) \\
\hline
\((9,11,21)\) & \((4800,440)=(3,110)\) & \(\{3,5\}\) & \(\mathbb{Z}[\sqrt{2},\sqrt{5}]\) \\
\hline
\((10,10,17)\) & \((3888,288)=(3,2)\) & \(\{2,3\}\) & \(\mathbb{Z}\) \\
\hline
\((10,12,12)\) & \((3267,143)=(11,143)\) & \(\{2,13\}\) & \(\mathbb{Z}[\sqrt{11}]\) \\
\hline
\((11,49,61)\) & \((72000,3720)=(5,930)\) & \(\{2,3\}\) & \(\mathbb{Z}[\sqrt{5},\sqrt{6}]\) \\
\hline
\((13,29,43)\) & \((35280,1848)=(5,462)\) & \(\{2,3,5,7\}\) & \(\mathbb{Z}[\sqrt{6},\sqrt{14}]\) \\
\hline
\((17,19,37)\) & \((25920,1368)=(5,38)\) & \(\{2,5\}\) & \(\mathbb{Z}[\sqrt{2}]\) \\
\hline
\((19,21,21)\) & \((18000,440)=(5,110)\) & \(\{2,5\}\) & \(\mathbb{Z}[\sqrt{10}]\) \\
\hline 
\caption{Classification of OI one-holed tori with associated invariants}
\label{tab:OIToriInvariants} 
\end{longtable}

\subsection{Constructing infinitely many commensurability classes of OI surfaces}\label{sec:ConstructingOIsurfaces}

In this section, we construct infinitely many pairwise non-commensurable OI surfaces. We start by making some observations. Consider the tiling of the hyperbolic plane $P(\mathcal{C}^+)$ by the orbit of $H(a,b,c)$ under the reflection groups $\Lambda$ generated by $\{ r_\alpha,r_\beta,r_\gamma,f_\alpha,f_\beta,f_\gamma\}$. Let $\mathcal{U}^\perp$ be the collection of geodesics in $P(\mathcal{C}^+)$ dual to orbits of ${\bf u}_*$ under $\Lambda$. By construction, distinct geodesics in $\mathcal{U}^\perp$ are disjoint. Their hyperbolic distances can be computed from the dual vector by the following formula
\[
\cosh(d_{\mathbb{H}^2}({\bf u}_1^\perp,{\bf u}_2^\perp)) = \frac{|B({\bf u}_1,{\bf u}_2)|}{\sqrt{B({\bf u}_1,{\bf u}_1)B({\bf u}_2,{\bf u}_2)}}.
\]
where ${\bf u}_i$ is the vector dual to ${\bf u}_i^\perp$. As the direct application of this formula, we have

\begin{lem}
\label{lem:ComputingCoshDistanceViaQuadraticForm}
Let ${\bf u} = (u_1,u_2,u_3)^{\top} \in \Lambda\cdot\{{\bf u}_\alpha,{\bf u}_\beta,{\bf u}_\gamma\}.$ Then we have
\[
\cosh(d_{\mathbb{H}^2}({\bf u}_\alpha^\perp,{\bf u}^\perp)) = u_1, \ \cosh(d_{\mathbb{H}^2}({\bf u}_\beta^\perp,{\bf u}^\perp)) = u_2, \  \cosh(d_{\mathbb{H}^2}({\bf u}_\gamma^\perp,{\bf u}^\perp)) = u_3 
\]
provided that ${\bf u}^\perp \neq {\bf u}_*^\perp$ in each case.
\end{lem}

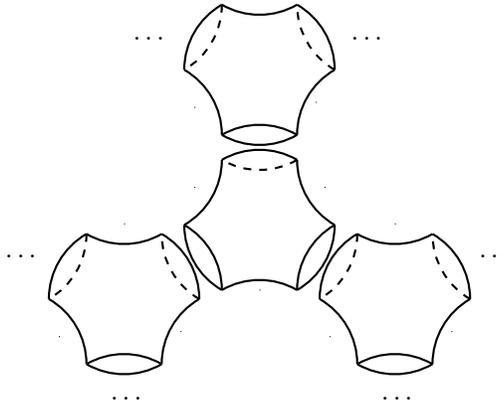
\begin{figure}[h!]
    \centering
    \begin{tikzpicture}[scale=1]
        % Original drawing at (0,0)
        \node (orig) at (0,0) {
            \begin{tikzpicture}[scale=1]
                \draw[thick] (1,0) arc (240:180:1);
                \draw[thick, dashed] (1/2,{sqrt(3)/2}) arc (300:240:1);
                \draw[thick] (-1/2,{sqrt(3)/2}) arc (0:-60:1);
                \draw[thick] (-1,0) arc (60:0:1);
                \draw[thick] (-1/2,{-sqrt(3)/2}) arc (120:60:1);
                \draw[thick] (1/2,{-sqrt(3)/2}) arc (180:120:1);
                %%%
                \draw[thick] (1/2,{sqrt(3)/2}) arc (60:120:1);
                \draw[thick] (-1,0) arc (180:240:1);
                \draw[thick] (1/2,{-sqrt(3)/2}) arc (300:360:1);
                %%%
                \draw ({sqrt(3)/2},1/2) -- ({sqrt(3)/2},1/2);
                \draw ({-sqrt(3)/2},1/2) -- ({-sqrt(3)/2},1/2);
                \draw (0,{-sqrt(3)/2}) -- (0,{-sqrt(3)/2});
            \end{tikzpicture}
        };
        % Second drawing: rotated by 45° and placed lower and closer to the first one.
        \node (rot) at (1.8,-1.1) {
            \begin{tikzpicture}[scale=1, rotate=60]
                \draw[thick] (1,0) arc (240:180:1);
                \draw[thick, dashed] (1/2,{sqrt(3)/2}) arc (300:240:1);
                \draw[thick] (-1/2,{sqrt(3)/2}) arc (0:-60:1);
                \draw[thick] (-1,0) arc (60:0:1);
                \draw[thick] (-1/2,{-sqrt(3)/2}) arc (120:60:1);
                \draw[thick, dashed] (1/2,{-sqrt(3)/2}) arc (180:120:1);
                %%%
                \draw[thick] (1/2,{sqrt(3)/2}) arc (60:120:1);
                \draw[thick] (-1,0) arc (180:240:1);
                \draw[thick] (1/2,{-sqrt(3)/2}) arc (300:360:1);
                %%%
                \draw ({sqrt(3)/2},1/2) -- ({sqrt(3)/2},1/2);
                \draw ({-sqrt(3)/2},1/2) -- ({-sqrt(3)/2},1/2);
                \draw (0,{-sqrt(3)/2}) -- (0,{-sqrt(3)/2});
            \end{tikzpicture}
        };

         % third drawing: rotated by 45° and placed lower and closer to the first one.
        \node (rot) at (0,1.95) {
            \begin{tikzpicture}[scale=1, rotate=60]
                \draw[thick] (1,0) arc (240:180:1);
                \draw[thick, dashed] (1/2,{sqrt(3)/2}) arc (300:240:1);
                \draw[thick] (-1/2,{sqrt(3)/2}) arc (0:-60:1);
                \draw[thick] (-1,0) arc (60:0:1);
                \draw[thick] (-1/2,{-sqrt(3)/2}) arc (120:60:1);
                \draw[thick, dashed] (1/2,{-sqrt(3)/2}) arc (180:120:1);
                %%%
                \draw[thick] (1/2,{sqrt(3)/2}) arc (60:120:1);
                \draw[thick] (-1,0) arc (180:240:1);
                \draw[thick] (1/2,{-sqrt(3)/2}) arc (300:360:1);
                %%%
                \draw ({sqrt(3)/2},1/2) -- ({sqrt(3)/2},1/2);
                \draw ({-sqrt(3)/2},1/2) -- ({-sqrt(3)/2},1/2);
                \draw (0,{-sqrt(3)/2}) -- (0,{-sqrt(3)/2});
            \end{tikzpicture}
        };

 % Fourth drawing: rotated by 45° and placed lower and closer to the first one.
        \node (rot) at (-1.8,-1.1) {
            \begin{tikzpicture}[scale=1, rotate=60]
                \draw[thick] (1,0) arc (240:180:1);
                \draw[thick, dashed] (1/2,{sqrt(3)/2}) arc (300:240:1);
                \draw[thick] (-1/2,{sqrt(3)/2}) arc (0:-60:1);
                \draw[thick] (-1,0) arc (60:0:1);
                \draw[thick] (-1/2,{-sqrt(3)/2}) arc (120:60:1);
                \draw[thick, dashed] (1/2,{-sqrt(3)/2}) arc (180:120:1);
                %%%
                \draw[thick] (1/2,{sqrt(3)/2}) arc (60:120:1);
                \draw[thick] (-1,0) arc (180:240:1);
                \draw[thick] (1/2,{-sqrt(3)/2}) arc (300:360:1);
                %%%
                \draw ({sqrt(3)/2},1/2) -- ({sqrt(3)/2},1/2);
                \draw ({-sqrt(3)/2},1/2) -- ({-sqrt(3)/2},1/2);
                \draw (0,{-sqrt(3)/2}) -- (0,{-sqrt(3)/2});
            \end{tikzpicture}
        };

        \node[anchor=west] at ($ (-1.8,2.4)$) {\ldots};
        
        \node[anchor=west] at ($(1.1,2.4)$) {\ldots};
        \node[anchor=west] at ($(1.5,- 2.4)$) {\ldots};
        \node[anchor=west] at ($(-2.1,- 2.4)$) {\ldots};
         \node[anchor=west] at ($(-3.5,- 0.5)$) {\ldots};
         \node[anchor=west] at ($ (2.8,- 0.5)$) {\ldots};
    \end{tikzpicture}
    \caption{Gluing copies of OI pair of pants $P(a,a,a)$ without twist, for a fixed $a\in \{2,3\}$.}
    \label{fig:GluingPants}
\end{figure}

The following proposition provides a source of infinitely many OI surfaces coming from gluing copies of $P(2,2,2)$ or $P(3,3,3)$, see \cref{fig:GluingPants}

\begin{prop}
\label{prop:InfinitelyManyOISurfaces}
    Let $P$ be either $P(2,2,2)$ or $P(3,3,3)$. Suppose $X$ is a surface obtained by gluing together finitely many copies of $P$ by reflections along the boundary components of the copies of $P$. Then $X$ is an ortho-integral surface.
\end{prop}

\begin{proof}
     By construction, the universal cover $\widetilde{X}$ of $X$ admits a tiling by an all right-angled hyperbolic hexagon $H(a,a,a)$ for $a \in \{2,3\}$. Thus, $\widetilde{X}$ can be identified with a convex subset of $P(\mathcal{C}^+)$ tilled by $H(a,a,a)$ for $a \in \{2,3\}$. The boundary of $\widetilde{X}$ consists of geodesics in $\mathcal{U}^\perp$. The set of cosh-lengths of orthogeodesics in $X$ is precisely the set of cosh-distances among geodesics in $\partial \widetilde{X}$. By \cref{lem:ComputingCoshDistanceViaQuadraticForm}, these cosh-distances are one of the coordinates of a vector in $\mathcal{U} = \Lambda\cdot \{{\bf u}_\alpha,{\bf u}_\beta,{\bf u}_\gamma\}$. For $a \in\{2,3\}$, the group $\Lambda$ has integral entries. That is, $\Lambda \leq \ORTH^+(q,\mathbb{Z})$. Therefore, $\mathcal{U} \subseteq \mathbb{Z}^3$. It follows that $X$ is an ortho-integral surface.  
\end{proof}

\begin{remark}\label{rem:NonEmptyOIForAnyType}
\cref{prop:InfinitelyManyOISurfaces} implies that the set of OI surfaces with fixed genus \(g\) and \(n>0\) boundary components is always nonempty, thereby completing the proof of \cref{thm:FinitenessOfOISurfaces}.
\end{remark}

Our examples of infinitely many non-commensurable OI surfaces come from the following family of surfaces. Fix $a \in \{2,3\}$. Define a sequence of hyperbolic surfaces $\{X_n\}$ of genus zero and $3\cdot 2^n$ boundary components inductively by letting $X_0 = P(a,a,a)$. The surface $X_{n+1}$ is obtained from $X_n$ by gluing one copy of $P(a,a,a)$ along each boundary component of $X_n$ by a reflection (see \cref{fig:GluingPants}).

It follows from \cref{prop:InfinitelyManyOISurfaces} that each \(X_n\) is an OI surface. Note that commensurable hyperbolic surfaces share an isometric finite-sheeted covering, so in the case of surfaces with geodesic boundary, they have isometric universal covers. We distinguish surfaces \(X_n\) up to commensurability by showing that their universal covers \(\widetilde{X}_n\) are pairwise non-isometric, as evidenced by differences in their inradii (i.e., the radii of the largest inscribed hyperbolic disks in each universal cover).

\begin{defi}
    Let $\Gamma$ be a finitely-generated non-elementary Fuchsian subgroup of $\PSL_2(\mathbb{R})$ such that the limit set $\Lambda_\Gamma \neq \partial_\infty \mathbb{H}^2$. Let $C(\Lambda_\Gamma)$ be the convex hull of $\Lambda_\Gamma$. The \emph{inradius} of $C(\Lambda_\Gamma)$ is defined to be
    \[
    r = \sup_{x \in C(\Lambda_\Gamma)}\inf_{\delta\in \partial C(\Lambda_\Gamma)}d_{\mathbb{H}^2}(x,\delta). 
    \]
\end{defi}
\noindent
Since $r(\Gamma)$ is an isometric invariant of the convex hull of the limit set of $\Gamma$, it is a commensurability invariant of $\Gamma$. 

\begin{lem}
If $\Gamma$ be a finitely-generated non-elementary Fuchsian subgroup of $\PSL_2(\mathbb{R})$ such that the limit set $\Lambda_\Gamma \neq \partial_\infty \mathbb{H}^2$, then $r(\Gamma)$ is finite.     
\end{lem}

\begin{proof}
    For a contradiction, suppose that there exists a sequence $\{x_n\} \subset C(\Lambda_\Gamma)$ such that \[\inf_{\delta\in \partial C(\Lambda_\Gamma)}d_{\mathbb{H}^2}(x_n,\delta) \to \infty.\] By assumption, $\Gamma$ admits a compact fundamental domain $D$ in $C(\Lambda_\Gamma)$ which is a polygon with finitely many sides. Since the function $x \mapsto \inf_{\delta\in \partial C(\Lambda_\Gamma)}d_{\mathbb{H}^2}(x,\delta)$ is $\Gamma$ invariant, we may assume that $\{x_n\} \subset D$. Passing to a subsequence, we suppose that $x_n \to x\in D$ since $D$ is compact. Let $\gamma$ be a boundary component in $\partial C(\Lambda_\Gamma)$ that realizes $\inf_{\delta\in \partial C(\Lambda_\Gamma)}d_{\mathbb{H}^2}(x,\delta)$. Then we have
    \[
    \inf_{\delta\in \partial C(\Lambda_\Gamma)}d_{\mathbb{H}^2}(x_n,\delta) \leq d_{\mathbb{H}^2}(x_n,\gamma) \leq d_{\mathbb{H}^2}(x_n,x) + d_{\mathbb{H}^2}(x,\gamma) + d_{\mathbb{H}^2}(\pi_\gamma(x),\pi_\gamma(x_n))
    \]
    where $\pi_\gamma: \mathbb{H}^2 \to \gamma$ is the orthogonal projection from $\mathbb{H}^2$ onto $\gamma$. As $n\to\infty$, the right-hand side tends to $d_{\mathbb{H}^2}(x,\gamma)$ while the left-hand side tends to $\infty$. This is the desired contradiction.
\end{proof}

\begin{lem}
\label{lem:NonIsometricUniCovers}
    Let $r_n$ be the inradius of $\widetilde{X}_n$. The sequence $r_n$ strictly increases and therefore $\widetilde{X}_n$ are non-isometric convex domain of $\mathbb{H}^2$. As a consequence, the surfaces $X_n$ are pairwise non-commensurable. 
\end{lem}

\begin{proof}
    It follows from the definition of $r_n$ and convexity of $\widetilde{X}_n$ that $r_n$ is the radius of the largest hyperbolic disk that can be embedded in $\widetilde{X}_n$. Let $x_n \in \widetilde{X}_n$ such that $\inf_{\delta\in \partial \widetilde{X}_n}d_{\mathbb{H}^2}(x_n,\delta)=r_n$. The universal cover $\widetilde{X}_{n+1}$ is obtained from $\widetilde{X}_n$ by gluing isometric copies of $\widetilde{X}_0$ along all boundary components of $\widetilde{X}_n$. Therefore in $\widetilde{X}_{n+1}$, there exists a disk of radius $r_n + \varepsilon$ centered at $x_n$ that remains embedded in $\widetilde{X}_{n+1}$ for a small enough $\varepsilon$. This shows that the inradius of $X_n$ is strictly increasing as $n \to \infty$.   
\end{proof}

\begin{proof}[Proof of \cref{thm:InfinitelyManyCommensurabilityClasses}]
    The sequence of hyperbolic surfaces $\{X_n\}$, defined as in the paragraph prior to \cref{fig:GluingPants}, consists of all OI surfaces by \cref{prop:InfinitelyManyOISurfaces}. The universal covers of these surfaces $X_n$ are pairwise non-isometric by \cref{lem:NonIsometricUniCovers}. It follows that the surfaces $\{X_n\}$ are pairwise non-commensurable.     
\end{proof}

\begin{remark}\label{rem:lastRemarkonCommensurability}
    For $a \in \{2,3\}$, the Fuchsian group $\Gamma_n$ that corresponds to the hyperbolic surface $X_n$ all contains $\Gamma_1$ up to conjugation. Since the invariant quaternion algebra of a non-elementary Fuchsian group is determined by two non-commuting hyperbolic elements, the groups $\Gamma_n$ all have the same invariant quaternion algebra. Among Fuchsian subgroups of finite co-area, the invariant quaternion algebra is the complete invariant of the commensurability class. It follows that the fact that $X_n$ are all pairwise non-commensurable is purely a metric phenomenon. 
\end{remark}

\begin{remark}
Thanks to the referee's suggestion, another way to show the sequence of OI surfaces $\{X_n\}$ are pairwise incommensurable is by comparing the orthokissing number and area:
$$
\operatorname{okiss}(X_n)=3\cdot 2^{n-1},\qquad
\operatorname{area}(X_n)=2\pi(3\cdot 2^{n}-2),
$$ for $n\ge 1$. If $X_n$ and $X_m$ are commensurable then the ratio $\operatorname{area}/\operatorname{okiss}$
coincides, hence
$$
\frac{3\cdot 2^{n}-2}{3\cdot 2^{n-1}}
=
\frac{3\cdot 2^{m}-2}{3\cdot 2^{m-1}}
\ \Rightarrow\ n=m.
$$
Thus, the sequence $\{X_n\}$ are pairwise incommensurable.
\end{remark}

Finally, as an application of \cref{prop:PantsAndOneHoledTorusFromRootFlipping} and \cref{lem:ComputingCoshDistanceViaQuadraticForm}, we construct examples of pairs of pants that are $d$-AOI but not OI using quadratic forms. 

\begin{prop}\label{prop:dAOI-not-OISurface}
The pair of pants $(\frac{n+1}{n},\frac{n+1}{n},\frac{n+1}{n})$ has the cosh-length spectrum lying in $\frac{1}{n}\mathbb{Z}$, where $n\in \mathbb{N}^{+}$. The pair of pants $(\frac{2n+3}{2n+1},\frac{2n+3}{2n+1},\frac{2n+3}{2n+1})$ has the cosh-length spectrum lying in $\frac{1}{2n+1}\mathbb{Z}$, where $n\in \mathbb{N}$.
\end{prop}

\begin{proof}
In the case of $P(\frac{n+1}{n},\frac{n+1}{n},\frac{n+1}{n})$, the group $G_P \leq \ORTH^+(q;\mathbb{Z})$ where $q$ is the quadratic form given in \cref{eq:QuadraticForm}. Furthermore, one can check that $G_P$ preserves the $\mathbb{Z}$-lattice $(\frac{1}{n}\mathbb{Z})^3$. This implies that the $G_P$-orbits of $\{{\bf u}_\alpha,{\bf u}_\beta,{\bf u}_\gamma\}$ are contained in $(\frac{1}{n}\mathbb{Z})^3$. Since the entries of the vectors in the $G_P$-orbits of $\{{\bf u}_\alpha,{\bf u}_\beta,{\bf u}_\gamma\}$ consist of the cosh-length spectrum of the  orthogeodesics of $P(\frac{n+1}{n},\frac{n+1}{n},\frac{n+1}{n})$. These pairs of pants provide examples of AOI surfaces but not OI surfaces. The argument is similar in the case of $P(\frac{2n+3}{2n+1},\frac{2n+3}{2n+1},\frac{2n+3}{2n+1})$.  
\end{proof}

\begin{remark}
Using the same argument as in \cref{prop:InfinitelyManyOISurfaces}, we can glue copies of the pair of pants $P(\frac{3}{2},\frac{3}{2},\frac{3}{2})$
together without any twists along their boundaries, either by gluing a copy to itself or to another, to obtain additional 2-AOI surfaces with higher topological complexity.
\end{remark}

%%%%%%%%%%%%%%%%%%%%%%%%%%%
% END of body
%%%%%%%%%%%%%%%%%%%%%%%%%%%
\bibliographystyle{alpha}
\bibliography{main}

\end{document}